\documentclass[12pt]{article}
\usepackage{fullpage}
\usepackage[english]{babel}
\usepackage[latin1]{inputenc}
\usepackage{epsfig}
\usepackage{color,graphicx,graphics,psfrag}
\usepackage{amsmath,amstext,amssymb,amsfonts, amscd}

\def\Box{\leavevmode\vbox{\hrule
     \hbox{\vrule\kern4pt\vbox{\kern4pt}%
           \vrule}\hrule}}
\def\blackbox{\leavevmode\vrule height 5pt width 4pt depth 0pt\relax}
\def\endproof{\null\hfill {$\blackbox$}\bigskip}

\newcounter{appendix}
\def\appendix{\advance\c@appendix by 1
   \def\thesection{\Alph{section}}
   \ifnum\c@appendix=1 \setcounter{section}{-1} \fi
   \@startsection {section}{1}{\z@}{-3.5ex plus -1ex minus 
   -.2ex}{2.3ex plus .2ex}{\Large\bf}}

\def\paragraph#1{{\bf #1\ }}

\newtheorem{lemma}{Lemma}

\newtheorem{definition}[lemma]{Definition}

\newtheorem{proposition}[lemma]{Proposition}

\newtheorem{remark}{Remark}[section]

\newcommand{\typsize}{L}

\title{Numerical approximation of the Euler-Poisson-Boltzmann model \\in the quasineutral limit} 
\author{P. Degond\footnote{Universit\'{e} de Toulouse; UPS, INSA, UT1, UTM and CNRS; 
Institut de Math\'{e}matiques de Toulouse;
F-31062 Toulouse, France. email: pierre.degond@math.univ-toulouse.fr}, H. Liu\footnote{Department of Mathematics, 
Iowa State University, Ames, IA 50011, USA. hliu@iastate.edu}, D. Savelief\footnote{Universit\'{e} de Toulouse; UPS, INSA, UT1, UTM and CNRS; 
Institut de Math\'{e}matiques de Toulouse;
F-31062 Toulouse, France. email: dominique.savelief@math.univ-toulouse.fr}, M-H. Vignal\footnote{Universit\'{e} de Toulouse; UPS, INSA, UT1, UTM and CNRS; 
Institut de Math\'{e}matiques de Toulouse;
F-31062 Toulouse, France. email: marie-helene.vignal@math.univ-toulouse.fr}} 
\date{} 
\begin{document}

\maketitle

\vspace{0.5 cm}
\begin{abstract}
This paper analyzes various schemes for the Euler-Poisson-Boltzmann (EPB) model of plasma physics. This model consists of the pressureless gas dynamics equations coupled with the Poisson equation and where the Boltzmann relation relates the potential to the electron density. If the quasi-neutral assumption is made, the Poisson equation is replaced by the constraint of zero local charge and the model reduces to the Isothermal Compressible Euler (ICE) model. We compare a numerical strategy based on the EPB model to a strategy using a reformulation (called REPB formulation). The REPB scheme captures the quasi-neutral limit more accurately. 
\end{abstract}

\medskip
\noindent
{\bf Acknowledgements:} P. D., D. S. and M-H. V. have been supported by the 'Fondation Sciences et Technologies pour 
l'Aéronautique et l'Espace', in the frame of the project 'Plasmax' (contract \# RTRA-STAE/2007/PF/002) and by the 
'F\'{e}d\'{e}ration de recherche sur la fusion par confinement magn\'{e}tique', 
in the frame of the contract 'APPLA' (Asymptotic-Preserving schemes for 
Plasma Transport) funded by the CEA (contract \# V3629.001 avenant 2). Liu's research was partially supported by the National Science Foundation under Kinetic FRG grant No. DMS 07-57227.


\vskip 0.4cm


\setcounter{equation}{0}
\section{Introduction}
\label{sec_intro}

The goal of this paper is to analyze various schemes for the Euler-Poisson-Boltzmann model of plasma physics. The Euler-Poisson-Boltzmann (EPB) model describes the plasma ions through a system of pressureless gas dynamics equations subjected to an electrostatic force. The electrostatic potential is related to the ion and electron densities through the Poisson equation. The Boltzmann relation provides a non-linear relationship between the potential and the electron density which allows to close the system. More precisely, the Euler-Poisson-Boltzmann model is written 
\begin{eqnarray}
& & \hspace{-1cm} \partial_t n + \nabla \cdot (nu) = 0, \label{EPB_n} \\
& & \hspace{-1cm} m ( \partial_t (nu) + \nabla \cdot (nu \otimes u)) = - e n \nabla \phi, \label{EPB_u} \\
& & \hspace{-1cm} - \Delta \phi = \frac{e}{{\epsilon}_0} (n - n^* \exp ( \frac{e \phi}{k_B T} )) . \label{EPB_phi}
\end{eqnarray}
Here, $n(x,t) \geq 0$, $u(x,t) \in {\mathbb R}^d $ and $\phi(x,t)\in {\mathbb R} $ stand for the ion density, ion velocity and electric potential respectively, which depend on the space-variable $x \in {\mathbb R}^d$ and on the time $t\geq0$. We suppose that the ions bear a single positive elementary charge $e$ and we denote by $m$, their mass. The electron temperature $T$ is supposed uniform and constant in time. $\epsilon_0$ and $k_B$ respectively refer to the vacuum permittivity and the Boltzmann constant. The operators $\nabla$, $\nabla \cdot $ and $\Delta$ are respectively the gradient, divergence and Laplace operators and $u \otimes u$ denotes the tensor product of the vector $u$ with itself. $n^*$ is usually fixed by either imposing zero total charge 
$$ \int_{{\mathbb R}^d} ( n(x,t) - n^* \exp ( \frac{e \phi(x,t)}{k_B T} )) \, dx = 0 , $$
or by assuming that the net charge is zero at a given point $x^*$ (for instance at the boundary): 
$$ n(x^*,t) - n^* \exp ( \frac{e \phi(x^*,t)}{k_B T} ) = 0 . $$
For simplicity, this work is restricted to dimension $d=1$ but the concepts extend to dimensions $d \geq 2$ without additional difficulties and numerical applications will be reported in future work.

The ion pressure force is neglected. This is a commonly made assumption in plasma physics \cite{Chen,Krall_Trivelpiece} for elementary text books of plasma physics. The inclusion of an ion pressure term would not modify the subsequent analysis and is omitted for simplicity. Additionally, the pressureless system has interesting multi-valued solutions which are lost in the case of the non-pressureless model \cite{Liu_Wang_1, Liu_Wang_2, Liu_Slemrod}. 

If the quasi-neutral assumption is made, the Poisson equation (\ref{EPB_phi}) is replaced by the constraint of zero local charge: 
$$ n = n^* \exp ( \frac{e \phi}{k_B T} ) , $$
In this context, we can write 
\begin{eqnarray}
& & \hspace{-1cm} n \nabla \phi = n^* \exp ( \frac{e \phi}{k_B T} ) \nabla \phi = \frac{k_B T}{e} \nabla ( n^* \exp ( \frac{e \phi}{k_B T} ) ) = \frac{k_B T}{e} \nabla n, \label{n_nabla_phi}
\end{eqnarray}
and the quasi-neutral Euler-Poisson-Boltzmann model coincides with Isothermal Compressible Euler (ICE) model: 
\begin{eqnarray*}
& & \hspace{-1cm} \partial_t n + \nabla \cdot (nu) = 0,  \\
& & \hspace{-1cm} m ( \partial_t (nu) + \nabla \cdot (nu \otimes u)) + \nabla (n k_B T) = 0 . \end{eqnarray*}

The passage from EPB to ICE can be understood by a suitable scaling of the model, which highlights the role of the scaled Debye length:
\begin{eqnarray}
& & \hspace{-1cm} \lambda = \frac{\lambda_D}{\typsize}, \quad \lambda_D = \left( \frac{\epsilon_0 k_B T}{e^2 n^*} \right)^{1/2}, \label{Debye} 
\end{eqnarray}
where $\typsize$ is the typical size of the system under consideration. $\lambda_D$ measures the spatial scale associated with the electrostatic interaction between the particles. The dimensionless parameter $\lambda$ is usually small, which formalizes the fact that the electrostatic interaction occurs at spatial scales which are much smaller than the usual scales of interest. However, there are situations, for instance in boundary layers, or at the plasma-vacuum interface, where the electrostatic interaction scale must be taken into account. This means that the choice of the relevant macroscopic length $\typsize$ may depend on the location inside the system and that in general, the parameter $\lambda$ may vary by orders of magnitude from one part of the domain to another one. The scaling will be presented in more detail in section \ref{sec_scaling}. 

This paper is concerned with discretization methods for the EPB model in situations where $\lambda$ can vary from order one to very small values. Therefore, the targeted schemes must correctly capture the transition from the EPB to the ICE models. With this objective in mind, we will compare two strategies: a first one which uses the EPB model in its original form, and a second one which reformulates the EPB model in such a way that it explicitly appears as a perturbation of the ICE model. This reformulation uses the Poisson equation in the form: 
$$  n = n^* \exp ( \frac{e \phi}{k_B T} ) -  \frac{{\epsilon}_0}{e} \Delta \phi. $$
Then, with the same algebra as for (\ref{n_nabla_phi}), we get:
\begin{eqnarray}
n \nabla \phi &=& n^* \exp ( \frac{e \phi}{k_B T} ) \nabla \phi -  \frac{{\epsilon}_0}{e} \Delta \phi \nabla \phi \nonumber \\
&=& \frac{k_B T}{e} \nabla ( n^* \exp ( \frac{e \phi}{k_B T} ) ) -  \frac{{\epsilon}_0}{e} (\nabla \cdot (\nabla \phi \otimes \nabla \phi) - \nabla ( \frac{ |\nabla \phi|^2}{2} ) ) \nonumber \\
& = & \frac{k_B T}{e} \nabla n + \frac{\epsilon_0 k_B T}{e^2} \nabla \Delta \phi -  \frac{{\epsilon}_0}{e} (\nabla \cdot (\nabla \phi \otimes \nabla \phi) - \nabla ( \frac{ |\nabla \phi|^2}{2} ) )  .  \label{n_nabla_phi_ref}
\end{eqnarray}
We note that this expression is reminiscent of the Maxwell stress tensor. Then, the EPB model is equivalently written as the following reformulated Euler-Poisson-Boltzmann (REPB) model: 
\begin{eqnarray}
& & \hspace{-1cm} \partial_t n + \nabla \cdot (nu) = 0, \label{REPB_n} \\
& & \hspace{-1cm} m ( \partial_t (nu) + \nabla \cdot (nu \otimes u)) + \nabla (n k_B T) = \nonumber \\
& & \hspace{2cm} = - {\epsilon}_0  \nabla ( \frac{k_B T}{e} \Delta \phi + \frac{ |\nabla \phi|^2}{2} ) +  {\epsilon}_0 \nabla \cdot (\nabla \phi \otimes \nabla \phi)   , \label{REPB_u} \\
& & \hspace{-1cm} - \Delta \phi = \frac{e}{{\epsilon}_0} (n - n^* \exp ( \frac{e \phi}{k_B T} )) . \label{REPB_phi}
\end{eqnarray}
In this way, the ICE appears at the left-hand side of (\ref{REPB_n}), (\ref{REPB_u}). The scaling analysis will show that the right-hand side of (\ref{REPB_u}) is of order $\lambda^2$ and is therefore negligible (if the gradients of the potential are smooth) in the limit $\lambda \to 0$. 

The goal of this paper is to propose and analyze two schemes for the EPB model which provide the correct ICE limit when $\lambda \to 0$: the first one is based on the initial formulation EPB and the second one, on the reformulated form REPB. To present the schemes, we note that both EPB and REPB can be put under the form 
\begin{eqnarray*}
& & \hspace{-1cm} \partial_t n + \nabla \cdot (nu) = 0,  \\
& & \hspace{-1cm} m ( \partial_t (nu) + \nabla \cdot (nu \otimes u)) + \nabla p(n) = S(n,\nabla \phi)   ,  \\
& & \hspace{-1cm} - \Delta \phi = \frac{e}{{\epsilon}_0} (n - n^* \exp ( \frac{e \phi}{k_B T} )) , 
\end{eqnarray*}
with 
$$ p(n) = 0, \quad S(n,\nabla \phi) = - e n \nabla \phi, $$
in the case of EPB and 
$$ p(n) = n k_B T, \quad S(n,\nabla \phi) = - {\epsilon}_0  \nabla ( \frac{k_B T}{e} \Delta \phi + \frac{ |\nabla \phi|^2}{2} ) +  {\epsilon}_0 \nabla \cdot (\nabla \phi \otimes \nabla \phi), $$
in the case of the REPB. 

Both schemes use the following time-semi-discretization which is implicit in the Poisson equation and in the source terms of the momentum equation:  
\begin{eqnarray*}
& & \hspace{-1cm} \delta^{-1} (n^{k+1} - n^k) + \nabla \cdot ((nu)^k) = 0,  \\
& & \hspace{-1cm} m ( \delta^{-1} ((nu)^{k+1} - (nu)^k) + \nabla \cdot ((nu)^k \otimes u^k)) + \nabla p(n^k) = S(n^{k+1},\nabla \phi^{k+1})   ,  \\
& & \hspace{-1cm} - \Delta \phi^{k+1} = \frac{e}{{\epsilon}_0} (n^{k+1} - n^* \exp ( \frac{e \phi^{k+1}}{k_B T} )) . 
\end{eqnarray*}
Here, $\delta$ is the time step and the exponent $k \in {\mathbb N}$ refers to the approximation at time $t^k = k \delta$ (i.e. $n^k(x) \approx n(x,t^k)$, \ldots). For the EPB model, this discretization is classical  \cite{Fabre_JCP_101_445}. For both systems, in spite of its implicit character, the recursion can be solved in an explicit way, by first updating the mass equation to find $n^{k+1}$, then using the Poisson equation to find the potential $\phi^{k+1}$ and finally using the momentum equation to find $u^{k+1}$. Of course, the space operators have to be discretized as well and a simple shock-capturing method is used, namely the Local Lax-Friedrichs or Rusanov scheme \cite{Leveque_2, Rusanov, Deg_Pey_Rus_Vil}. A time-splitting is used where the hydrodynamic equations are evolved without the source terms. Then, the Poisson equation is solved with the value of the density found at the end of the first step of the splitting. With the newly computed value of the potential, the evolution of the hydrodynamic quantities due to the source terms are computed. 

In section \ref{sec_time}, we will show that both schemes are Asymptotic-Preserving. The Asympto\-tic-Preserving property can be defined as follows. Consider a singular perturbation problem $P^\lambda$ whose solutions converge to those of a limit problem $P^0$ when $\lambda \to 0$ (here $P^\lambda$ is the EPB model and $P^0$ is the ICE model). A scheme $P^\lambda_{\delta,h}$ for problem $P^\lambda$ with time-step $\delta$ and space-step $h$ is called Asymptotic Preserving (or AP) if it is stable independently of the value of $\lambda$ when $\lambda \to 0$ and if the scheme $P^0_{\delta,h}$ obtained by letting $\lambda \to 0$ in $P^\lambda_{\delta,h}$ with fixed $(\delta,h)$ is consistent with problem $P^0$. This property is illustrated by the commutative diagram  below: 
$$ \begin{CD} 
P^\lambda_{\delta,h} @>{(\delta,h) \to 0}>>  P^\lambda \\
@VV{\lambda \to 0}V  @VV{\lambda \to 0}V\\
P^0_{\delta,h} @>{(\delta,h) \to 0}>> P^0 
\end{CD}
$$
The concept of an AP scheme has been introduced by S. Jin \cite{Jin} for diffusive limits of kinetic models and has been widely expanded since then. 

In section \ref{sec_time}, we will perform a linearized stability analysis which shows that both schemes are stable independently of $\lambda$ in the limit $\delta \to 0$, provided that the usual CFL condition of the ICE model is satisfied. However, the scheme based on the REPB formulation has several advantages  over the one based on the EPB form. A first advantage lies in the fact that the hydrodynamic part of the EPB model is a pressureless gas dynamics model, which is weakly unstable and can produce delta concentrations (see e.g. \cite{Bouchut, Bouchut_James, Brenier_Grenier, CL1, CL2}). By contrast, the hydrodynamic part of the REPB is the usual ICE model, which is strongly hyperbolic, and which is much more stable than the pressureless gas dynamics model. A second advantage is the fact that the limit $\lambda \to 0$ of the REPB-based scheme provides a conservative discretization of the ICE model, while that of the EPB-based scheme leads to a scheme in non-conservative form. We can expect that the accuracy of the PB-based scheme degrades when $\lambda \ll 1$ for solutions involving discontinuities, which is not the case for the REPB-based scheme.

To illustrate the theoretical findings, one-dimensional numerical simulations are presented in section \ref{sec_num}, following a discussion of the spatial discretization in section \ref{sec_spatial}. 
First, setting $\lambda = 1$, analytic solutions can be derived in the form of solitary waves thanks to the Sagdeev potential theory \cite{Chen}. Both schemes are compared to these analytical solutions, and show a similar behavior, with a slightly larger numerical diffusion in the case of the REPB scheme. Then a Riemann problem test case consisting of two outgoing shock waves is investigated. In the $ \lambda \ll 1 $ regime, the REPB scheme captures the right hydrodynamic shocks while the EPB scheme develops spurious oscillations. This confirms the better behavior of the REPB scheme in the $ \lambda \ll 1 $ regime. Finally, a test-problem related to multivalued solutions and proposed in \cite{Liu_Wang_1,Liu_Wang_2} is investigated. In this case, both schemes shows a similar behavior. To summarize, the REPB scheme captures well the $\lambda \ll 1$ limit but is slightly more diffusive in the $ \lambda = O(1) $ regime. But the extra numerical diffusion is mild and then shows that the REPB scheme is superior to and should be preferred over the EPB scheme in most situations.

We conclude this section by some bibliographical remarks. The Euler-Poisson-Boltzmann model has been recently analyzed in the context of sheath dynamics \cite{Liu_Slemrod} and multi-valued solutions have been computed using level set methods \cite{Liu_Wang_1, Liu_Wang_2}. Numerical schemes for the quasineutral limit of plasma problems has been the subject of vast literature, mostly for the Vlasov-Poisson equation and for Particle-in-Cell (PIC) methods. It is virtually impossible to cite all relevant references and we only refer to the seminal ones \cite{Cohen_JCP_46_15, Langdon_JCP_51_107, Mason_JCP_41_233, Mason_JCP_51_484}. For fluid models of plasmas, the literature is comparatively less abundant. We can refer to the pioneering work \cite{Fabre_JCP_101_445}, and more recently to \cite{Choe_JCP_170_550, Collela_JCP_149_168, Schneider_IJNM_05_399, Shumlak_JCP_187_620}. Recently, AP-schemes for the two-fluid Euler-Poisson model \cite{Cri_Deg_Vig_07, Deg_Liu_Vig_08} or Vlasov-Poisson model
\cite{BCDS_09, DDNSV_10} in the quasineutral limit have been proposed. The drift-fluid limit of magnetized plasmas has also been considered \cite{Deg_Del_Neg_10, DDSV_09} as well as other applications such as small Mach-number flows \cite{Deg_Tan_10}. However, none of these works is concerned with Boltzmannian electrons.

\setcounter{equation}{0}
\section{The EPB model and the scaling}
\label{sec_scaling}

In this section we present a derivation of the EPB model from the two fluid Euler-Poisson
system and we introduce its scaling. From now on, we will restrict ourselves to one-dimensional models.

\subsection{Derivation of the EPB model}
\label{subsec_derivation}

We consider a  plasma composed of two species of charged particles: positively charged ions and electrons. The ions are supposed singly charged. The modeling of such a plasma by means of fluid equations uses a system of compressible Euler equations for each species, coupled by the Poisson equation. 
We denote by $ m_{i,e} $ the ion and electron masses, 
 $ n_{i,e} $ the ion and electron densities and
 $ u_{i,e} $ the ion and electrons mean velocities.
We denote by $ e $ the elementary charge (i.e. the ion charge is $+e>0$
and the electron charge is $-e<0$).
We assume that the ion temperature is negligible so that
the ions  follow a pressureless gas dynamics model.
Electrons are assumed isothermal with a
 non-zero constant and uniform temperature $ T_e $.
Then the electron pressure law satisfies
 $ p_e = n_e k_B T_e $
where $ k_B $ denotes the Boltzmann constant.
 The balance laws for both species are given by
\begin{eqnarray}
 \label{epbif-dim-iden}
& & \hspace{-1cm}    \partial_t n_i + 
   \partial_x ( n_i u_i )  = 0 , \\
 \label{epbif-dim-imom}
& & \hspace{-1cm}     m_i (\partial_t ( n_i u_i ) +
   \partial_x ( n_i u_i^2 ))  = 
   - e n_i \partial_x \phi , \\
 \label{epbif-dim-eden}
& & \hspace{-1cm}   \partial_t n_e + 
   \partial_x ( n_e u_e )  = 0 , \\
 \label{epbif-dim-emom}
& & \hspace{-1cm}  m_e (\partial_t ( n_e u_e ) +
   \partial_x ( n_e u_e^2 )) +  \partial_x ( n_e k_B T_e )  = 
    e n_e \partial_x \phi ,
\end{eqnarray}
where $ \phi $ denotes the electric potential. It satisfies
 the Poisson equation:
\begin{equation}
 \label{bif-dim-poisson}
  - \epsilon_0 \partial_{x}^2 \phi = 
   e ( n_i - n_e ),
\end{equation}
where $ \epsilon_0 $ is the vacuum permittivity.

The electrons being much lighter than the ions, it is legitimate to take the limit $ m_e \to 0 $ in the electron momentum equation. In this limit, we formally obtain:
$$  \partial_x ( n_e k_B T_e ) = 
   e n_e \partial_x \phi.
$$
Integration with respect to $ x $ leads to the Boltzmann relation:
\begin{equation}
 \label{boltz-dim}
  n_e = n^* \exp \left( \frac{ e \phi }{ k_B T_e } \right),
\end{equation}
where $ n^* $ is fixed by some condition (e.g. vanishing total charge or vanishing local charge at one given point such as a boundary point, see section \ref{sec_intro}). The Boltzmann relation shows that the electron density automatically adjusts
 to the potential. The EPB model therefore consists of the ion mass and momentum conservation equations (\ref{epbif-dim-iden}), (\ref{epbif-dim-imom}), the Poisson equation (\ref{bif-dim-poisson}) and the Boltzmann relation (\ref{boltz-dim}). With a change of notation $n_i \to n$, $u_i \to u$, $m_i \to m$, $T_e \to T$, we find the EPB model (\ref{EPB_n})-(\ref{EPB_phi}) which has been introduced in section \ref{sec_intro}. 
 
We note that, in the present one-dimensional setting, the electron velocity $u_e$ can be computed from the electron density equation (\ref{epbif-dim-eden}), thanks to the value of $n_e$ and therefore, of $\phi$, obtained by the resolution of the EPB model. Therefore, the computation of $u_e$ is decoupled from the computation of the other unknowns $n_i$, $u_i$ and $\phi$ involved in the EPB model and will be discarded in the present work. In two or higher dimensions, the computation of $u_e$ requires the resolution of the electron momentum equation, which takes the form in the small electron mass limit: 
\begin{equation}
\label{emom_small_mass}
  \partial_t ( n_e u_e ) +
   \nabla \cdot ( n_e u_e \otimes u_e ) +  n_e \nabla \psi  = 0, 
\end{equation}
where $\psi = \lim_{m_e \to 0} ( m_e^{-1} (k_B T \ln n_e - e \phi)) $. The quantity $\psi$ plays the same role as the pressure in the incompressible Euler equation. It is computed thanks to the electron density equation (\ref{epbif-dim-eden}) which appears as a (non-zero) divergence constraint on $u_e$. In this sense, the limit $m_e \to 0$ is similar to the 'low Mach-number' limit of isentropic compressible gas dynamics. However, the question of the resolution of (\ref{emom_small_mass}) is left to future work.

\subsection{Scaling of the EPB model}
\label{subsec_scaling}

In this section, we return to the EPB model in the form (\ref{EPB_n})-(\ref{EPB_phi}) and, with the notations of section \ref{sec_intro}, we introduce a scaling of the physical quantities. Let $x_0$, $t_0$, $u_0$, $\phi_0$ and $n_0$ be space, time, velocity, potential and density scales. Scaled position, time, velocity, potential and density are defined by $ \bar{x} = x/x_0$, $ \bar{t} = t/t_0$, $ \bar{u} = u/u_0$ , $ \bar{\phi} = - \phi/\phi_0$ and  $ \bar{n} = n/n_0$. We choose $x_0$ to be the typical size of the system (for instance an inter-electrode distance or the size of the vacuum chamber). The velocity scale is chosen equal to the ion sound speed 
$ u_0 = ( k_B T / m )^{1/2} $. We note that the ion sound speed is constructed with the ion mass but with the electron temperature. This is clear from the ICE model (see section \ref{sec_intro}). We also choose $n_0=n^*$. Finally, $\phi_0 = k_B T /e$ is the so-called thermal potential. Note that we have introduced a sign change in the potential scaling because we find it more convenient to work in terms of the electron potential energy rather than in terms of the electric potential. 

Inserting this scaling and omitting the bars gives rise to the EPB model in scaled form: 
\begin{eqnarray}
 \label{nop-den}
& & \hspace{-1cm}  \partial_t n^\lambda + \partial_x ( n^\lambda u^\lambda )  = 0, \\
 \label{nop-mom}
& & \hspace{-1cm}   \partial_t ( n^\lambda u^\lambda ) + 
  \partial_x ( n^\lambda u^\lambda u^\lambda )  = 
    n^\lambda \partial_x \phi^\lambda, \\
 \label{nop-pot}
& & \hspace{-1cm}  \lambda^2 \partial_{x}^2 \phi^\lambda  = n^\lambda - e^{-\phi^\lambda}.
\end{eqnarray}
where $\lambda$ is the scaled Debye length (\ref{Debye}). 

It will be useful to consider the linearized EPB model about the state defined by $n^\lambda = 1$, $u^\lambda = 0$, $\phi^\lambda = 0$ (which is obviously a stationary solution). Expanding $n^\lambda = 1 + \varepsilon \tilde n^\lambda$, $u^\lambda = \varepsilon \tilde u^\lambda$ and $\phi^\lambda = \varepsilon \tilde \phi^\lambda$, with $\varepsilon \ll 1$ being the intensity of the perturbation to the stationary state, and retaining only the linear terms in $\varepsilon$, we find the linearized EPB model: 
\begin{eqnarray}
& & \hspace{-1cm} \partial_t \tilde n^\lambda + \partial_x \tilde u^\lambda = 0, \label{LEPB_n} \\
& & \hspace{-1cm} \partial_t \tilde u^\lambda = \partial_x \tilde \phi^\lambda , \label{LEPB_u} \\
& & \hspace{-1cm} \lambda ^2 \partial^2_x \tilde \phi^\lambda = \tilde n^\lambda + \tilde \phi^\lambda. \label{LEPB_phi} 
\end{eqnarray}
Introducing $\hat n^\lambda$, $\hat u^\lambda$, $\hat \phi^\lambda$, the partial Fourier transforms of $\tilde n^\lambda$, $\tilde u^\lambda$, $\tilde \phi^\lambda$ with respect to $x$, we are led to the system of ODE's: 
\begin{eqnarray}
& & \hspace{-1cm} \partial_t \hat n^\lambda + i \xi \hat u^\lambda = 0, \label{FLEPB_n} \\
& & \hspace{-1cm} \partial_t \hat u^\lambda = i \xi \hat \phi^\lambda , \label{FLEPB_u} \\
& & \hspace{-1cm} - \lambda ^2 \xi^2 \hat \phi^\lambda = \hat n^\lambda + \hat \phi^\lambda , \label{FLEPB_phi} 
\end{eqnarray}
where $\xi$ is the Fourier dual variable to $x$. We note that the general solution of this model takes the form
\begin{eqnarray}
& & \hspace{-1cm} \left( \begin{array}{c} \hat n^\lambda \\ \hat u^\lambda \end{array} \right) = \sum_{\pm} e^{s^\lambda_\pm t} \left( \begin{array}{c} \hat n^\lambda_\pm \\ \hat u^\lambda_\pm \end{array} \right), 
\label{EPB_gen_sol}
\end{eqnarray}
with
\begin{eqnarray}
& & \hspace{-1cm} s^\lambda_\pm  = \pm \frac{i \xi}{(1 + \lambda^2 \xi^2)^{1/2}}, 
\label{EPB_s}
\end{eqnarray}
and $n^\lambda_\pm$, $u^\lambda_\pm$ are given functions of $\xi$, fixed by the initial conditions of the problem. In particular, since $s^\lambda_\pm $ are pure imaginary numbers, the $L^2$ norm of the solution is preserved with time.

Now, we investigate the quasi-neutral limit $\lambda \to 0$ in the next section.

\subsection{The quasineutral limit: the ICE model}
\label{subsec_quasineutral}

Formally passing to the limit $ \lambda \to 0 $ in the EPB model in scaled form and supposing that $n^\lambda \to n^0$, $u^\lambda \to u^0$, $\phi^\lambda \to \phi^0$, we are led to the following model: 
\begin{eqnarray*}
& & \hspace{-1cm} \partial_t n^0 + \partial_x ( n^0 u^0 )  = 0, \\
& & \hspace{-1cm}  \partial_t ( n^0 u^0 ) + \partial_x ( n^0 u^0 u^0 )  = 
n^0 \partial_x \phi^0, \\
& & \hspace{-1cm}  0 = n^0 - e^{-\phi^0}.
\end{eqnarray*}
As a consequence of the last relation (which imposes to the ions to satisfy the Boltzmann relation of the electrons), we can write (see also (\ref{n_nabla_phi})):
\begin{equation}
 n^0 \partial_x \phi^0 = 
  e^{-\phi^0} \partial_x \phi^0 = 
  - \partial_x \left( e^{-\phi^0} \right) = 
  - \partial_x n^0,
  \label{n_nabla_phi_1D}
\end{equation}
and, inserting this relation into the momentum equation leads to: 
$$ \partial_t ( n^0 u^0 ) +
  \partial_x ( n^0 u^0 u^0 ) +
 \partial_x n^0 = 0.
$$
Therefore, the quasineutral limit $\lambda \to 0$ consists of the Isothermal Compressible Euler system (ICE) complemented by the Boltzmann relation for the potential:
\begin{eqnarray*}
& & \hspace{-1cm}    \partial_t n^0 + \partial_x ( n^0 u^0 )  = 0, \\
& & \hspace{-1cm}    \partial_t ( n^0 u^0 ) + 
    \partial_x ( n^0 u^0 u^0 ) +
    \partial_x n^0  = 0, \\
& & \hspace{-1cm}    n^0  = e^{-\phi^0}.
\end{eqnarray*}

Similarly to the EPB model, the ICE model can be linearized about the state defined by $n^0 = 1$, $u^0 = 0$, $\phi^0 = 0$. We find (with the same notations as for the EPB model), in Fourier space: 
\begin{eqnarray*}
& & \hspace{-1cm} \partial_t \hat n^0 + i \xi \hat u^0 = 0,  \\
& & \hspace{-1cm} \partial_t \hat u^0 + i \xi \hat n^0 = 0 .  
\end{eqnarray*}
The general solution of this model takes the same form (\ref{EPB_gen_sol}) as for the linearized EPB model with $s^\lambda_\pm$ replaced by $ s^0_\pm  = \pm i \xi. $ We note that (see (\ref{EPB_s}))
$$ s^0_\pm = \lim_{\lambda \to 0} s^\lambda_\pm. $$ 
Therefore, the wave speeds of the linearized EPB model converge to those of the linearized ICE model (which are nothing but the acoustic wave speeds). There is no singularity of the limit $\lambda \to 0$ as regards the wave-speeds. In this respect the quasineutral limit $\lambda \to 0$ is not a singular limit. This fact contrasts with the situation of the quasineutral limit of the 2-fluid Euler system, where the electron plasma oscillation frequency converges to infinity. Therefore, we expect that the numerical treatment of the quasineutral limit in the Euler-Poisson-Boltzmann case will be easier.

\subsection{Reformulation of the EPB model}
\label{subsec_reformulation}

To better capture the transition from the EPB model to the ICE model, it is useful to reformulate the EPB model in such a way that it explicitly appears as a perturbation of the ICE model. Using (\ref{nop-pot}), in the spirit of (\ref{n_nabla_phi_1D}), we can write (see also (\ref{n_nabla_phi_ref})): 
\begin{eqnarray*}
 n^\lambda \partial_x \phi^\lambda &=&
  \left( e^{-\phi^\lambda}+\lambda^2 \partial_{xx}^2 \phi^\lambda \right)
   \partial_x \phi^\lambda  \\
  &=& \partial_x \left( - e^{-\phi^\lambda} + \frac{\lambda^2}{2} ( \partial_x \phi^\lambda )^2 \right) \\
  &=& \partial_x
  \left( \lambda^2 \partial_{xx}^2 \phi^\lambda - 
  n^\lambda + \frac{\lambda^2}{2} ( \partial_x \phi^\lambda )^2\right).
\end{eqnarray*}
We note that some simplification arises in the 1-dimensional case, compared to the multi-dimensional case of (\ref{n_nabla_phi_ref}). Inserting this expression in the momentum equation leads to the reformulated EPB systems (REPB): 
\begin{eqnarray}
 \label{ref-den}
& & \hspace{-1cm}  \partial_t n^\lambda + \partial_x ( n^\lambda u^\lambda )  = 0, \\
 \label{ref-mom}
& & \hspace{-1cm}  \partial_t ( n^\lambda u^\lambda ) + 
   \partial_x ( n^\lambda u^\lambda u^\lambda ) +
   \partial_x n^\lambda  = \lambda^2 \partial_x
   \left( \partial_{xx}^2 \phi^\lambda +
   \frac{1}{2} ( \partial_x \phi^\lambda )^2 \right), \\
 \label{ref-pot}
& & \hspace{-1cm}  n^\lambda - e^{-\phi^\lambda}  = \lambda^2 \partial_{x}^2 \phi^\lambda .
\end{eqnarray}

In this formulation, the ICE model explicitly appears at the left-hand side of (\ref{ref-den})-(\ref{ref-pot}).  Additionally, the remaining terms, at the right-hand side of the equations are formally of order $\lambda^2$. Therefore, the EPB model explicitly appears as an order $O(\lambda^2)$ perturbation of the ICE model. 

We stress the fact that the REPB model is {\bf equivalent} to the original EPB model. However, at the discrete level, schemes based on the REPB model may differ from those based on the EPB model. The goal of the present article is to compare the properties of schemes based on these two formulations, in relation to their Asymptotic-Preserving (AP) properties when $\lambda \to 0$. 

\begin{remark}
A formal expansion (using a Chapman-Enskog methodology) up to the first order in $ \lambda^2 $ of the EPB model leads to the following model:
\begin{eqnarray*}
 & & \hspace{-1cm} \partial_{t} n^{\lambda} + \partial_{x} (n^{\lambda} u^{\lambda} ) = 0, \\
 & & \hspace{-1cm} \partial_{t} (n^{\lambda} u^{\lambda} ) + \partial_{x} ( n^{\lambda} u^{\lambda} u^{\lambda} )
 + \partial_{x} n^{\lambda} = - \lambda^2 n^{\lambda} \partial_{x} 
 \left( \frac{1}{n^{\lambda}} \partial_{x}^{2} ( \log n^{\lambda} ) \right) + O(\lambda^4).
\end{eqnarray*}
We see that the EPB is a perturbation of the ICE model by a dispersive term with a third order derivative in $n^{\lambda}$. For this reason, in the sequel, the $ \lambda = O(1) $ regime will be referred to as the dispersive regime, while the $ \lambda \ll 1 $ will be referred to as the hydrodynamic regime.
\end{remark}

\setcounter{equation}{0}
\section{Time semi-discretization, AP property and linearized stability}
\label{sec_time}

\subsection{Time-semi-discretization and AP property}
\label{subsec_time_disc}

We denote by $ \delta  $ the time step. For any function $g(x,t)$, we denote by
$ g^m (x) $ an approximation of $g(x,t^m)$  with $ t^m = m \delta $. We present two time-semi-discretizations of the problem. The first one is based on the EPB formulation, and the second one, on the REPB formulation.

\subsubsection{Time-semi-discretization based on the EPB formulation}
\label{subsubsec_time_EPB}

Classically, when dealing with Euler-Poisson problems, the force term in the momentum equation is taken implicitly. In the case of the two-fluid model (when the electrons are modeled by the compressible Euler equations instead of being described by the Boltzmann relation), S. Fabre has shown that this implicitness is needed for the stability of the scheme (an explicit treatment of the force term leads to an unconditionally unstable scheme \cite{Fabre_JCP_101_445}). Additionally, this implicitness still gives rise to an explicit resolution, since the mass conservation can be used to update the density, then the Poisson equation is used to update the potential, and finally the resulting potential is inserted in the momentum equation to update the velocity. 

We will reproduce this strategy here and consider the following time-semi-discretization based on the EPB formulation: 
\begin{eqnarray*}
& & \hspace{-1cm}  \delta^{-1} (n^{\lambda, m+1} - n^{\lambda, m}) + 
  \partial_x ( n^{\lambda, m} u^{\lambda, m} ) = 0, \\
& & \hspace{-1cm}   \delta^{-1} ( n^{\lambda, m+1} u^{\lambda, m+1} - n^{\lambda, m} u^{\lambda, m})
  + \partial_x ( n^{\lambda, m} u^{\lambda, m} u^{\lambda, m} )  =
  n^{\lambda, m+1} \partial_x \phi^{\lambda, m+1}, \\
& & \hspace{-1cm}   \lambda^2 \partial_{x}^2 \phi^{\lambda, m+1}  =
  n^{\lambda, m+1} - e^{-\phi^{\lambda, m+1}}.
\end{eqnarray*}

This scheme is Asymptotic-Preserving. Indeed, letting $\lambda \to 0$ in this scheme with a fixed $\delta$ leads to 
\begin{eqnarray}
 \nonumber
& & \hspace{-1cm}   \delta^{-1} (n^{0, m+1} - n^{0, m}) + 
  \partial_x ( n^{0, m} u^{0, m} ) = 0, \\
 \nonumber
& & \hspace{-1cm}   \delta^{-1} ( n^{0, m+1} u^{0, m+1} - n^{0, m} u^{0, m})
  + \partial_x ( n^{0, m} u^{0, m} u^{0, m} )  =
  n^{0, m+1} \partial_x \phi^{0, m+1}, \\
 \label{class-pot-0}
& & \hspace{-1cm}   0  =
  n^{0, m+1} - e^{-\phi^{0, m+1}}.
\end{eqnarray}
By using (\ref{class-pot-0}) and the same algebra as for (\ref{n_nabla_phi_1D}), we find that this scheme is equivalent to 
\begin{eqnarray*}
& & \hspace{-1cm}   \delta^{-1} (n^{0, m+1} - n^{0, m}) + 
  \partial_x ( n^{0, m} u^{0, m} ) = 0, \\
& & \hspace{-1cm}   \delta^{-1} ( n^{0, m+1} u^{0, m+1} - n^{0, m} u^{0, m})
  + \partial_x ( n^{0, m} u^{0, m} u^{0, m} ) + \partial_x ( n^{0, m+1} ) =
  0, \\
& & \hspace{-1cm}   0  =
  n^{0, m+1} - e^{-\phi^{0, m+1}},
\end{eqnarray*}
which provides a semi-implicit discretization of the ICE model, with an implicit treatment of the pressure term in the momentum conservation equation. 

However, when the scheme is discretized in space, the algebra leading to (\ref{n_nabla_phi_1D}) is no longer exact. Let us denote by $D \phi^{0, m+1}$ the discretization of the space derivative operator $\partial_x \phi^{0, m+1}$. Then, the limit $\lambda \to 0$ of the fully discrete scheme gives rise to the approximation $n^{0, m+1}  D (\ln n^{0, m+1})$ of the space derivative $\partial_x ( n^{0, m+1} )$ instead of the natural derivative $D  n^{0, m+1}$. In particular, this expression is not in conservative form. Therefore, the use of this scheme may lead to a wrong shock speed if shock waves are present in the solution.

Another drawback of this scheme is the lack of pressure term in the momentum equation. As a consequence, the hydrodynamic part of the model is a pressureless gas dynamics model, which is a weakly ill-posed model (with, e.g. the possibility of forming delta concentrations \cite{Bouchut, Bouchut_James, Brenier_Grenier, CL1, CL2}). The weak instability may lead to spurious oscillations in the solution. 

For these reasons, another scheme, based on the REPB formulation, is proposed in the next section.

\subsubsection{Time-semi-discretization based on the REPB formulation}
\label{subsubsec_time_REPB}

We reproduce the same strategy (i.e. an implicit evaluation of the force term in the momentum conservation equation) starting from the REPB formulation. This leads to the following scheme: 
\begin{eqnarray}
 \label{refs-den} 
& & \hspace{-1cm}    \delta^{-1} (n^{\lambda, m+1} - n^{\lambda, m})
   + \partial_x ( n^{\lambda, m} u^{\lambda, m} )  = 0, \\
  \label{refs-mom}
& & \hspace{-1cm}     \delta^{-1} (( n^{\lambda, m+1} u^{\lambda, m+1} )
    - ( n^{\lambda, m} u^{\lambda, m} ))
    + \partial_x ( n^{\lambda, m} u^{\lambda, m} u^{\lambda, m})  + \partial_x n^{\lambda, m}  = \nonumber \\
& & \hspace{6cm}     =  \lambda^2 \partial_x \left( \partial_{x}^2 \phi^{\lambda, m+1} + 
    \frac{1}{2} ( \partial_x \phi^{\lambda, m+1} )^2 \right), \\
 \label{refs-pot}
& & \hspace{-1cm}    \lambda^2 \partial_{x}^2 \phi^{\lambda, m+1}  = 
  n^{\lambda, m+1} - e^{-\phi^{\lambda, m+1}}.
\end{eqnarray}

Formally passing to the limit $ \lambda \to 0 $ with fixed $\delta$
in this scheme leads to the following scheme: 
\begin{eqnarray}
 \label{refs-den-0} 
& & \hspace{-1cm}    \delta^{-1} (n^{0, m+1} - n^{0, m})
   + \partial_x ( n^{0, m} u^{0, m} )  = 0, \\
  \label{refs-mom-0}
& & \hspace{-1cm}     \delta^{-1} (( n^{0, m+1} u^{0, m+1} )
    - ( n^{0, m} u^{0, m} ))
    + \partial_x ( n^{0, m} u^{0, m} u^{0, m})  + \partial_x n^{0, m}  = 0, \\
 \label{refs-pot-0}
& & \hspace{-1cm}    0  = 
  n^{0, m+1} - e^{-\phi^{0, m+1}}.
\end{eqnarray}
Eqs. (\ref{refs-den-0}), (\ref{refs-mom-0}) are the standard time-semi-discretization of the ICE model. We now note that the pressure term $\partial_x n^{0, m}$ is explicit (it was implicit in the scheme based on the EPB formulation). Additionally, if a space discretization is used, the discretization of this term will stay in conservative form, by contrast to the EPB-based scheme. Finally, the hydrodynamic part of the scheme (\ref{refs-den})-(\ref{refs-pot}) is based on the ICE model, not on the pressureless gas dynamics model. Therefore, its discretization will avoid the possible spurious oscillations that might appear in the EPB-based scheme in the presence of discontinuities or sharp gradients.

\subsection{Linearized stability analysis of the time-semi-discretization}
\label{sub_sec_stability}

The goal of this section is to analyze the linearized stability properties of both schemes. More precisely, we want to show that both schemes are stable under the CFL condition of the ICE model, irrespective of the value of $\lambda$ when $\lambda \to 0$. This property is known as 'Asymptotic-Stability' and is a component of the Asymptotic-Preserving property (see section \ref{sec_intro}). Indeed, the faculty of letting $\lambda \to 0$ in the scheme with fixed $\delta$ is possible only if the stability condition of the scheme is independent of $\lambda$ in this limit. We will prove $L^2$-stability uniformly with respect to $\lambda$ for the linearization of the problem (\ref{FLEPB_n})-(\ref{FLEPB_phi}).  

In general, time semi-discretizations of hyperbolic problems are unconditionally unstable. This is easily verified on the discretization (\ref{refs-den-0})-(\ref{refs-pot-0}) of the ICE model. This is because the skew adjoint operator $\partial_x$ has the same effect as a centered space-differencing. For fully discrete schemes, stability is obtained at the price of adding numerical viscosity. To mimic the effect of this viscosity, in the present section, we will consider the linearized Viscous Euler-Poisson-Boltzmann (VEPB) model, which consists of the linearized EPB model (\ref{LEPB_n})-(\ref{LEPB_phi}) with additional viscosity terms (in this section, we drop the tildes for notational convenience): 
\begin{eqnarray*}
& & \hspace{-1cm} \partial_t  n^\lambda + \partial_x  u^\lambda - \beta \partial_x^2 n^\lambda = 0, \\
& & \hspace{-1cm} \partial_t  u^\lambda - \beta \partial_x^2 u^\lambda = \partial_x  \phi^\lambda , \\
& & \hspace{-1cm} \lambda ^2 \partial^2_x  \phi^\lambda =  n^\lambda +  \phi^\lambda.
\end{eqnarray*}
where $\beta$ is a numerical viscosity coefficient. We keep in mind that, in the spatially discretized case, $\beta$ proportional to the mesh size $h$: 
\begin{equation}
\beta = c h , \label{beta}
\end{equation}
with the constant $c$ to be specified later on. Similarly, the linearized Reformulated Viscous Euler-Poisson-Boltzmann (RVEPB) model is written: 
\begin{eqnarray*}
& & \hspace{-1cm} \partial_t  n^\lambda + \partial_x  u^\lambda - \beta \partial_x^2 n^\lambda = 0,   \\
& & \hspace{-1cm} \partial_t  u^\lambda + \partial_x  n^\lambda - \beta \partial_x^2 u^\lambda = \lambda^2 \partial_x^3  \phi^\lambda ,  \\
& & \hspace{-1cm} \lambda ^2 \partial^2_x  \phi^\lambda =  n^\lambda +  \phi^\lambda. 
\end{eqnarray*}

The time discretization of these two formulations (which are also linearizations of the EPB or REPB-based schemes with added viscosity terms) are given by 
\begin{eqnarray}
& & \hspace{-1cm} \delta^{-1}   (n^{\lambda,m+1} - n^{\lambda,m})  + \partial_x  u^{\lambda,m} - \beta \partial_x^2 n^{\lambda,m} = 0, \label{DLVEPB_n} \\
& & \hspace{-1cm} \delta^{-1}   (u^{\lambda,m+1} - u^{\lambda,m}) - \beta \partial_x^2 u^{\lambda,m} = \partial_x  \phi^{\lambda,m+1} ,\label{DLVEPB_u} \\
& & \hspace{-1cm} \lambda ^2 \partial^2_x  \phi^{\lambda,m+1} =  n^{\lambda,m+1} +  \phi^{\lambda,m+1}, \label{DLVEPB_phi}
\end{eqnarray}
for the EPB-based scheme and by 
\begin{eqnarray}
& & \hspace{-1cm} \delta^{-1}   (n^{\lambda,m+1} - n^{\lambda,m}) + \partial_x  u^{\lambda,m} - \beta \partial_x^2 n^{\lambda,m} = 0, \label{DLRVEPB_n} \\
& & \hspace{-1cm} \delta^{-1}   (u^{\lambda,m+1} - u^{\lambda,m}) + \partial_x  n^{\lambda,m} - \beta \partial_x^2 u^{\lambda,m} = \lambda^2 \partial_x^3  \phi^{\lambda,m+1} , \label{DLRVEPB_u} \\
& & \hspace{-1cm} \lambda ^2 \partial^2_x  \phi^{\lambda,m+1} =  n^{\lambda,m+1} +  \phi^{\lambda,m+1}, \label{DLRVEPB_phi} 
\end{eqnarray}
for the REPB-based one. 

Passing to Fourier space with $\xi$ being the dual variable to $x$, and eliminating $\hat \phi^{\lambda,m+1}$, we find the following recursion relations:  
\begin{eqnarray*}
& & \hspace{-1cm} \delta^{-1}   (\hat n^{\lambda,m+1} - \hat n^{\lambda,m})  + i \xi \hat  u^{\lambda,m} + \beta \xi^2 \hat n^{\lambda,m} = 0,  \\
& & \hspace{-1cm} \delta^{-1}   (\hat u^{\lambda,m+1} -  \hat u^{\lambda,m}) + \frac{i \xi}{1+\lambda^2 \xi^2} \hat n^{\lambda,m+1} + \beta \xi^2 \hat u^{\lambda,m} = 0 , 
\end{eqnarray*}
for the EPB-based scheme and 
\begin{eqnarray*}
& & \hspace{-1cm} \delta^{-1}   (\hat n^{\lambda,m+1} - \hat n^{\lambda,m}) + i \xi \hat u^{\lambda,m} + \beta \xi^2 \hat n^{\lambda,m} = 0,  \\
& & \hspace{-1cm} \delta^{-1}   ( \hat u^{\lambda,m+1} -  \hat u^{\lambda,m}) + i \xi \hat  n^{\lambda,m} - \frac{i \lambda^2 \xi^3}{1+\lambda^2 \xi^2} \hat n^{\lambda,m+1} + \beta \xi^2 \hat u^{\lambda,m} = 0 , 
\end{eqnarray*}
for the REPB-based one.

The characteristic equations for these two recursion formulas are
\begin{eqnarray}
& & \hspace{-1cm} 
q^2 - 2q (1 - \beta \xi^2 \delta - \frac{\xi^2 \delta^2}{2 (1+ \lambda^2 \xi^2)} ) + (1 - \beta \xi^2 \delta)^2 = 0
, \label{char_EPB} 
\end{eqnarray}
and 
\begin{eqnarray}
& & \hspace{-1cm} 
q^2 - 2q (1 - \beta \xi^2 \delta + \frac{\lambda^2 \xi^4 \delta^2}{2 (1+ \lambda^2 \xi^2)} ) + (1 - \beta \xi^2 \delta)^2 + \xi^2 \delta^2  = 0
, \label{char_REPB}
\end{eqnarray}
respectively, where $q$ is the characteristic root. Each of these quadratic equations has two roots $q^\lambda_\pm(\xi)$ which provide the two independent solutions of the corresponding recursion formulas. Their most general solution is of the form
\begin{eqnarray*}
& & \hspace{-1cm} \left( \begin{array}{c} \hat n^{\lambda,m}(\xi) \\ \hat u^{\lambda,m}(\xi) \end{array} \right) = \sum_{\pm} (q^\lambda_\pm(\xi))^m \left( \begin{array}{c} \hat n^\lambda_\pm(\xi) \\ \hat u^\lambda_\pm(\xi) \end{array} \right), \quad \forall m \in {\mathbb N}, 
\end{eqnarray*}
where $n^\lambda_\pm(\xi)$ and $u^\lambda_\pm(\xi)$ depend on the initial condition only. A necessary and sufficient condition for $L^2$ stability is that $|q^\lambda_\pm(\xi)| <1$. However, requesting this condition for all $\xi \in {\mathbb R}$ is too restrictive. To account for the effect of a spatial discretization in this analysis, we must restrict the range of admissible Fourier wave-vectors $\xi$ to the interval $[-\frac{\pi}{h}, \frac{\pi}{h}]$. Indeed, a space discretization of step $h$ cannot represent wave-vectors of magnitude larger than $\frac{\pi}{h}$. This motivates the following definition of stability: 

\begin{definition}
The scheme is stable if and only if
\begin{equation}
|q^\lambda_\pm(\xi)|\leq 1, \quad \forall \xi \quad \mbox{ such that } \quad |\xi| < \frac{\pi}{h}.
\label{stab_cnd}
\end{equation}
\label{def_stab}
\end{definition}

Now, our goal is to find sufficient conditions on $\delta$ such that either schemes are stable. More precisely, we prove: 

\begin{proposition}
For both the EPB-based scheme (\ref{DLVEPB_n})-(\ref{DLVEPB_phi}) or the REPB-based scheme (\ref{DLRVEPB_n})-(\ref{DLRVEPB_phi}), there exists a constant $C>0$ independent of $\lambda$ when $\lambda \to 0$ such that if $\delta \leq C h$, the scheme is stable. 
\label{prop_stab}
\end{proposition}

This condition states that the schemes are stable irrespective of how small $\lambda$ is. We say that the schemes are 'Asymptotically-Stable' in the limit $\lambda \to 0$. We note that this stability condition is similar to the CFL condition of the ICE model, which is the limit model when $\lambda \to 0$.   

\medskip
\noindent
{\bf Proof:} We first define conditions such that the constant term of the quadratic equations (\ref{char_EPB}) or (\ref{char_REPB}) is less than $1$. For the EPB-based scheme (see \ref{char_EPB}), this condition is 
$\delta \leq {1}/{\beta \xi^2}$, for all $\xi$ such that $|\xi| \leq {\pi}/{h}$. For reasons which will become clear below, we rather impose:
\begin{eqnarray}
& & \hspace{-1cm} 
\delta \leq \frac{1}{2 \beta \xi^2}, \quad \forall \xi \quad \mbox{ such that } \quad |\xi| \leq \frac{\pi}{h}  
, \label{cond_EPB} 
\end{eqnarray}
which, with (\ref{beta}), is equivalent to:  
\begin{eqnarray}
& & \hspace{-1cm} 
\delta \leq C_1 h, \quad \mbox{ with } \quad C_1 = \frac{1}{2 c \pi^2}
. \label{cond2_EPB} 
\end{eqnarray}
For the REPB-based scheme (see \ref{char_REPB}), this condition is 
$$\delta \leq \frac{2 \beta}{\beta^2 \xi^2 + 1}, \quad \forall \xi \quad \mbox{ such that } \quad |\xi| \leq \frac{\pi}{h}  
,  $$
or, with (\ref{beta}):  
\begin{eqnarray}
& & \hspace{-1cm} 
\delta \leq  C_1 h, \quad \mbox{ with } \quad C_1 = \frac{2c}{1+c^2 \pi^2}
. \label{cond2_REPB} 
\end{eqnarray}

Now, under these conditions, we are guaranteed that the two roots satisfy (\ref{stab_cnd}) if and only if the discriminant of the quadratic equation is non-positive. Indeed, in this condition, the two roots are conjugate complex numbers and their product, i.e. the square of their module, which is equal to the constant term of the quadratic equation, is less than one. 

It is a matter of computation to check that the discriminant has the same sign as the expression $F(\delta)$ given by:
\begin{eqnarray} 
& & \hspace{-1cm} 
F(\delta) = \delta^2 + 4 \beta (1 + \lambda^2 \xi^2) \delta - 4 \frac{1 + \lambda^2 \xi^2}{\xi^2} 
, \label{Delta_EPB} 
\end{eqnarray}
in the case of the EPB-based scheme and by 
$$F(\delta) = \delta^2 - 4 \beta \frac{1 + \lambda^2 \xi^2}{\lambda^2 \xi^2} \delta - 4 \frac{1 + \lambda^2 \xi^2}{\lambda^4 \xi^6}
, $$
in the case of the REPB-based scheme. 

In the case of the EPB-based scheme, we use (\ref{cond_EPB}) to estimate the second term of $F(\delta)$ in (\ref{Delta_EPB}): 
\begin{eqnarray*} 
& & \hspace{-1cm} 
F(\delta) \leq \delta^2  - 2 \frac{1 + \lambda^2 \xi^2}{\xi^2} , 
\end{eqnarray*}
and a sufficient condition for $F(\delta)$ to be non-positive is that $\delta \leq  \sqrt 2 {(1 + \lambda^2 \xi^2)^{1/2}}{|\xi|}^{-1}$. This relation is true  for all $\xi$  such that $ |\xi| \leq {\pi}/{h}$ if  
$\delta \leq \sqrt 2 \pi^{-1} (h^2 + \lambda^2 \pi^2)^{1/2}$ and a sufficient condition is that $\delta \leq C_2 h$, with $C_2=\sqrt 2 \pi^{-1}$.  Now, taking $C = \min\{C_1,C_2\}$ with $C_1$ given by (\ref{cond2_EPB}) leads to the result. In fact, the optimal numerical viscosity is such that $C_1=C_2$, i.e. $c= (2 \sqrt 2 \pi)^{-1}$. 
In the case of the REPB scheme, we estimate $F(\delta)$ by 
\begin{eqnarray*} 
& & \hspace{-1cm} 
F(\delta) \leq \delta^2 - 4 \beta \frac{1 + \lambda^2 \xi^2}{\lambda^2 \xi^2} \delta 
,  
\end{eqnarray*}
and a sufficient condition for $F(\delta)$ to be non-positive is that $\delta \leq  4 \beta {(1 + \lambda^2 \xi^2)}({\lambda^2 \xi^2})^{-1}$. In view of (\ref{beta}), this relation is true  for all $\xi$  such that $ |\xi| \leq {\pi}/{h}$ if  
$\delta \leq 4 c (\lambda^2 \pi^2)^{-1} h (h^2 + \lambda^2 \pi^2)^{1/2}$ and a sufficient condition is that $\delta \leq C_2 h$, with $C_2=4c$.  Now, taking $C = \min\{C_1,C_2\}= C_1$ with $C_1$ given by (\ref{cond2_REPB}) leads to the result. The optimal numerical viscosity can be chosen to minimize $C$, which leads to  $c= \pi^{-1}$. This ends the proof of statement \ref{prop_stab}. \endproof

As a conclusion, we can see that both the EPB-based and REPB-based schemes have similar Asymptotic-Stability properties as $\lambda \to 0$. Therefore, they must be selected on the basis of other criteria. The fact that the REPB-based scheme has a well-posed hydrodynamic part and leads to a discretization of the ICE model in conservative form in the limit $\lambda \to 0$ are indications that this scheme should be preferred to the EPB-based scheme. In the next section, we will present numerical results that support this statement.

\setcounter{equation}{0}
\section{Spatial discretization}
\label{sec_spatial}

 We introduce $ ( C_j )_{j=1}^{N} $ a uniform subdivision of
 the computational domain $ \Omega \in \mathbb{R} $  such that 
 $ \Omega = \cup_{j=1}^{N} C_j $.
 The interface between $ C_{j} $ and $ C_{j+1} $ is the point $ x_{j+1/2} $.
We denote by $ U_{j}^{m} $ the approximate 
 vector of the density and momentum at time~$ t^{m} $ on the cell~$ C_j $,
$$ U_{j}^{m} = \begin{pmatrix}
  n_{j}^{m} \\
  (n u)_{j}^{m}
 \end{pmatrix}.
$$
We use a time-splitting method to compute the
 density and momentum at time $ t^{m+1} $.
 There are three steps to pass from  $ U^{m} $ to $ U^{m+1} $ which are described below.

\subsection{Hydrodynamic part}
\label{subsec_hydro}

 The first step of the splitting is the finite-volume computation of the state
 $ U^{\#} $ such that
$$ \frac{ U^{\#} - U_{j}^{m} }{ \delta } +
 \frac{ F_{j+1/2}^{m} - F_{j-1/2}^{m} }{ h } =  0,
$$
where $ F_{j+1/2} $ is the numerical flux computed
 at the interface $ x_{j+1/2} $.

We have used a Local Lax-Friedrichs \cite{Leveque_2} (or Rusanov \cite{Rusanov} or degree 0 polynomial \cite{Deg_Pey_Rus_Vil}) solver. This solver is an improved version of the Lax-Friedrich solver
 which has been successfully used in conjunction with AP-schemes for the two-fluid Euler-Poisson problem (see \cite{Cri_Deg_Vig_07})  or for small Mach-number flows \cite{Deg_Tan_10}. This solver depends on a local estimate of the maximal characteristic speed. This estimate proceeds as follows. We introduce
$$ (a^{+})_{j+1/2}^{m} = \max \left( u_{j+1/2}^m + 1 , u_{j+1}^m + 1 \right)  \quad \mbox{ and } \quad  (a^{-})_{j+1/2}^{m} = \min \left( u_j^m - 1 , u_{j+1/2}^m - 1 \right) , $$
where $ u_{j+1/2}^m = ( u_j^m + u_{j+1}^m )/2 $.
Then, the local maximal characteristic velocity is estimated by
$$ a_{j+1/2}^m =  \max \left( | (a^{-})_{j+1/2}^{m} | ,    | (a^{+})_{j+1/2}^{m} | \right),  $$ 
and the numerical flux at $ x_{j+1/2} $ is given by:
\begin{equation}
 F_{j+1/2} = \frac{1}{2} \left( F(U_j^m) + F(U_{j+1}^m)
 + a_{j+1/2}^m \left( U_j^m - U_{j+1}^m \right) \right).
\label{num_flux}
\end{equation}

The time step $ \delta $ must satisfy the CFL condition
$ \frac{ \delta }{ h } \leq \max a_{j+1/2}^m $
to ensure stability.
In practice, the time step is chosen at each iteration
to enforce this stability condition. As for boundary conditions, we impose fictitious states $U_l$ and $U_r$ across the left and right boundaries respectively and compute the corresponding fluxes across the boundaries using the same formula 
(\ref{num_flux}).

\subsection{Potential update}
\label{subsec_pot}

The second step of the time-splitting is the computation of the potential.
 We use $ n^{\#} $ to compute $ \phi^{m+1} $ with the discretized
 Poisson equation given by a finite difference approximation of the Poisson equation.
$$ \lambda^2 h^{-2}  (\phi_{j-1}^{m+1} - 2 \phi_{j}^{m+1}
                                               + \phi_{j+1}^{m+1})
  + e^{ - \phi_{j}^{m+1} } = n_{j}^{\#}.
$$
 The boundary conditions are given by
 $ \phi_l = - \log n_{l}^{\#} $ on
 the left hand side of the domain
 and $ \phi_r = - \log n_{r}^{\#} $
 on the right hand side of the domain.
The non-linear system is solved with newton method.
This iterative algorithm needs a good initial guess of the solution to
 be efficient.
The initial guess is the potential at previous time $ \phi^{m} $.
For the first step, we choose the quasi-neutral potential
 $ ( \phi_{j}^{0} ) = ( - \log n_{j}^0 ) $ as an initial guess.

\subsection{Source term}
\label{subsec_source}

In the REPB form, the source term $\mathcal{Q}$ at the right-hand side of (\ref{refs-mom}) can be written according to:
$$ \mathcal{Q} =  \lambda^2 \left( \partial_x ( \partial_{xx}^2 \phi)  +  
 \partial_x \phi \,  \partial^2_x \phi \right) , $$
This expression is discretized thanks to a finite difference approximation.
For a cell $ C_{j} $ in the domain, the source term $ \mathcal{Q}_j $ is given by
centered finite difference approximation:
\begin{eqnarray*}
\mathcal{Q}_{j} = \frac{\lambda^2}{2h^3} \left(  \left(
 \phi_{j+2} - 2 \phi_{j+1} + 2 \phi_{j-1} - \phi_{j-2} \right)
 + \left( \phi_{j+1} - 2 \phi_{j} + \phi_{j-1} \right)
   \left( \phi_{j+1} - \phi_{j-1} \right) \right).
\end{eqnarray*}
On the first cell $ C_{1} $ the source term is computed
 using a decentered finite difference approximation:
\begin{eqnarray*}
\mathcal{Q}_{1} = \frac{\lambda^2}{h^3} \left( \left(
 \phi_{3} - 3 \phi_{2} + 3 \phi_{1} - \phi_{l} \right)
 + \frac{1}{2} \left( \phi_{2} - 2 \phi_{1} + \phi_{l} \right)   \left( \phi_{2} - \phi_{l} \right) \right),
\end{eqnarray*}
and similarly in the last cell.


\setcounter{equation}{0}
\section{Numerical results}
\label{sec_num}

We present three  classes of numerical results :
 the first test case is a solitary wave travelling in a plasma.
This test case shows the ability of the numerical schemes to handle the dispersive regime.
The second test case is related to the quasi-neutral limit $ \lambda \to 0 $ of the Euler-Poisson-Boltzmann model:
 it is a Riemann problem to check the ability of the scheme to handle hydrodynamic phenomena like shocks.
The last test case has been previously investigated by Liu and Wang in \cite{Liu_Wang_1,Liu_Wang_2} and corresponds to the occurrence of multi-valued solutions in the semi-classical setting.



\subsection{Soliton test case}
\label{subsec_soliton}

\subsubsection{Description}
\label{subsubsec_soliton_description}

The dispersive nature of the Euler-Poisson system is shown in \cite{Liu_Slemrod}. 
Therefore, the EPB system, like other nonlinear dispersive models such as the KdV equation, exhibits solitary
wave solutions. These special solutions are particularly convenient 
to test the ability of the EPB and REPB schemes
to capture the dispersive regime. Solitary waves also provide interesting quantitative checks. Indeed, while travelling through the plasma, the soliton maintains its shape and velocity.
Therefore, one can check the accuracy of the numerical schemes by observing how well 
they preserve the soliton shape and velocity over time.

We now summarize the establishment of this special solution. 
We refer to \cite{Chen} for a detailed description and 
physical considerations. For this derivation, we use non-dimensional units and we now precise the corresponding scaling units.
The space scale related to these solitons is the Debye length $ \lambda_{D} $.
For this reason, we take $ \lambda = 1 $ in all this section.
The size of the computational domain is equal to several
 Debye lengths (about $ 50 \lambda_{D} $ are used in the subsequent simulations).
The density of the undisturbed plasma (away from the support of the solitary wave) is
 chosen as the characteristic density and in dimensionless units, is equal to $ 1 $, so that the electron density is equal to $e^{-\phi}$ where $\phi$ is the electrostatic potential energy.

In the frame moving with the wave, we denote by  $ n_{s}, u_{s}, \phi_{s} $
 the density, velocity and potential of the plasma.
These quantities are constant in time, and satisfy the following
 relations:
\begin{eqnarray*}
 &  & \partial_x ( n_s u_s )  = 0, \\
 &  & \partial_x ( n_s u_s u_s )  = n_s \partial_x \phi_s, \\
 &  & \partial_{x}^{2} \phi_s  = n_s - e^{-\phi_s}.
\end{eqnarray*}
The momentum being uniform in $ x $, we write $ q = n_s u_s = n_0 u_0$ where $n_0 = n_s(0)$ and $u_0 = u_s(0)$. 
The momentum conservation law can be written as $ \partial_x q^2 / n_s = n_s \partial_x \phi_s $.
Consequently, we have $ \partial_x \left( {q^2}/{n_s^2} \right) = 2 \partial_x \phi_s $.
For all $ x \in [0,x_{\max}] $, we get:
\begin{equation}
 \frac{1}{2} \left( \frac{q^2}{n_s^2} \right)(x) -
 \frac{1}{2} \left( \frac{q^2}{n_s^2} \right)(0) =
 \phi_s ( x ) - \phi_s ( 0 ).
\end{equation}
The potential being defined up to an additive constant, we choose this constant such that
$ \phi_s (0) = 0 $.  
The ion density is then given by 
\begin{equation}
\label{eq:density_sagdeev}
 n_s (x) = \left( \frac{1}{n_0^2} + \frac{2 \phi_s}{ n_0^2 u_0^2 } \right)^{-1/2},
\end{equation}
In the present analysis, we assume that $ n_{0} = 1 $, i.e. $n_0$ is equal to the density of the undisturbed plasma.
Inserting this relation in the Poisson
 equation yields the following equation
 for the potential:
\begin{equation} \label{sol_poisson_eq}
 \partial_{x}^{2} \phi_s
 = \left( 1 + \frac{2 \phi_s}{ u_0^2 } \right)^{-1/2} - e^{-\phi_s},
\end{equation}
with the condition
\begin{equation} \label{sol_poisson_eq_x=0}
\phi_s (0) = 0 .
\end{equation}
One needs another condition at $ x = 0 $ to set up a Cauchy problem for (\ref{sol_poisson_eq}).
Note that $ \phi_s \equiv 0 $ is
 an obvious solution of equation (\ref{sol_poisson_eq}) and satisfies the homogeneous condition $\partial_x \phi_s ( 0 ) = 0 $.
It corresponds to the state of the undisturbed plasma.
To capture a non-trivial solution, we must
 alter this condition by a small
 disturbance, setting it to 
\begin{equation} \label{sol_poisson_eq_x=1}
 \partial_x \phi_s ( 0 ) = \eta  , 
\end{equation}
with $\eta$ 'small'.
 
The behavior of these solutions 
 can be clarified thanks to the theory of the Sagdeev potential
 (which is a primitive with respect to $\phi_s$ of the right-hand side of (\ref{sol_poisson_eq})).
Details on this study can be found in \cite{Chen}.
Here we just recall that
shock waves in a cold-ion plasma can exist
 only for $ 1 < u_0 < 1.6 $ (Bohm criterion).
The sign of $ \eta $ determines
 the type of solution which can be found:
 a positive $ \eta $ leads to potential barrier that forms a sheath, 
 whereas a negative value gives rise to a monotonic transition to a negative $ \phi $ which forms a solitary wave corresponding to a potential
 and density disturbance propagating to the
 right (for instance) with velocity $ u_0 $.
There is no analytic solution to equation (\ref{sol_poisson_eq}).
Resorting to numerical simulation is the only way to determine these solutions.
Details about this numerical method are given below.
Once the soliton potential is known,
 the density is computed thanks to (\ref{eq:density_sagdeev}).
Since the present analysis has been performed in a co-moving frame with the soliton, moving with velocity $u_0$, the velocity in the laboratory frame is $u_s + u_0$ with $u_s = q/n_s$. 
In the subsequent numerical simulations $ n_{s}, u_{0}+u_{s}, \phi_{s} $ are taken as
 an initial condition for the scheme.

\subsubsection{Numerical results for the soliton test case}
\label{subsubsec_results_soliton}

In this test case, where the Debye length and the size of the computational domain are of the same order, both the EPB and REPB schemes are expected to be correct.
However, this test case provides a way to achieve quantitative comparisons of the numerical solutions to an analytical reference solution. These comparisons permit to quantify the order of accuracy and amount of numerical diffusion of the two schemes.
The analytical solution is easily obtained at time $ t $ by
 a spatial translation of the solution at time $0$ of a distance $ u_0 t $.
 
This test case is implemented as follows. 
The boundary conditions are periodic, which ensures that the
 soliton can be tracked on long simulation times without
the need for a large computational domain.
The length $L$ of the computational domain is defined in relation to the choice of the initial condition as explained below.
We denote by $ t_{L} = L / u_{0} $ the travel time of the soliton
in the domain $ [0,L] $ .
First a numerical solution of (\ref{sol_poisson_eq}) with initial conditions (\ref{sol_poisson_eq_x=0}), (\ref{sol_poisson_eq_x=1}) is computed
with an explicit finite difference scheme on a mesh of step $ \Delta {x}^{\text{ref}} $ and provides the reference solution.
$ \Delta {x}^{\text{ref}} $ is chosen small enough to provide an almost 'exact' reference solution. 

For suitable $ \eta $ and $ u_{0} $ given by the study of the Sagdeev potential \cite{Chen},
 the potential oscillates in space for positive $ x $.
One wants a single potential well to initialize the scheme.
To this aim the number of nodes $ N^{\text{ref}} $ is taken such that the initial condition shows a single peak.
Moreover, it is such that $ \phi^{\text{ref}}_{N^{\text{ref}}} $ is close
enough to $ 0 $ to ensure that periodic boundary conditions will be accurate enough.
This reference solution is interpolated to provide an initial condition
 for the EPB and REPB schemes, and to compute numerical errors on the density,
 momentum and potential.
 
\begin{figure}[hbtp]
 \begin{minipage}[c]{.46\linewidth}
 \includegraphics[scale=0.4]{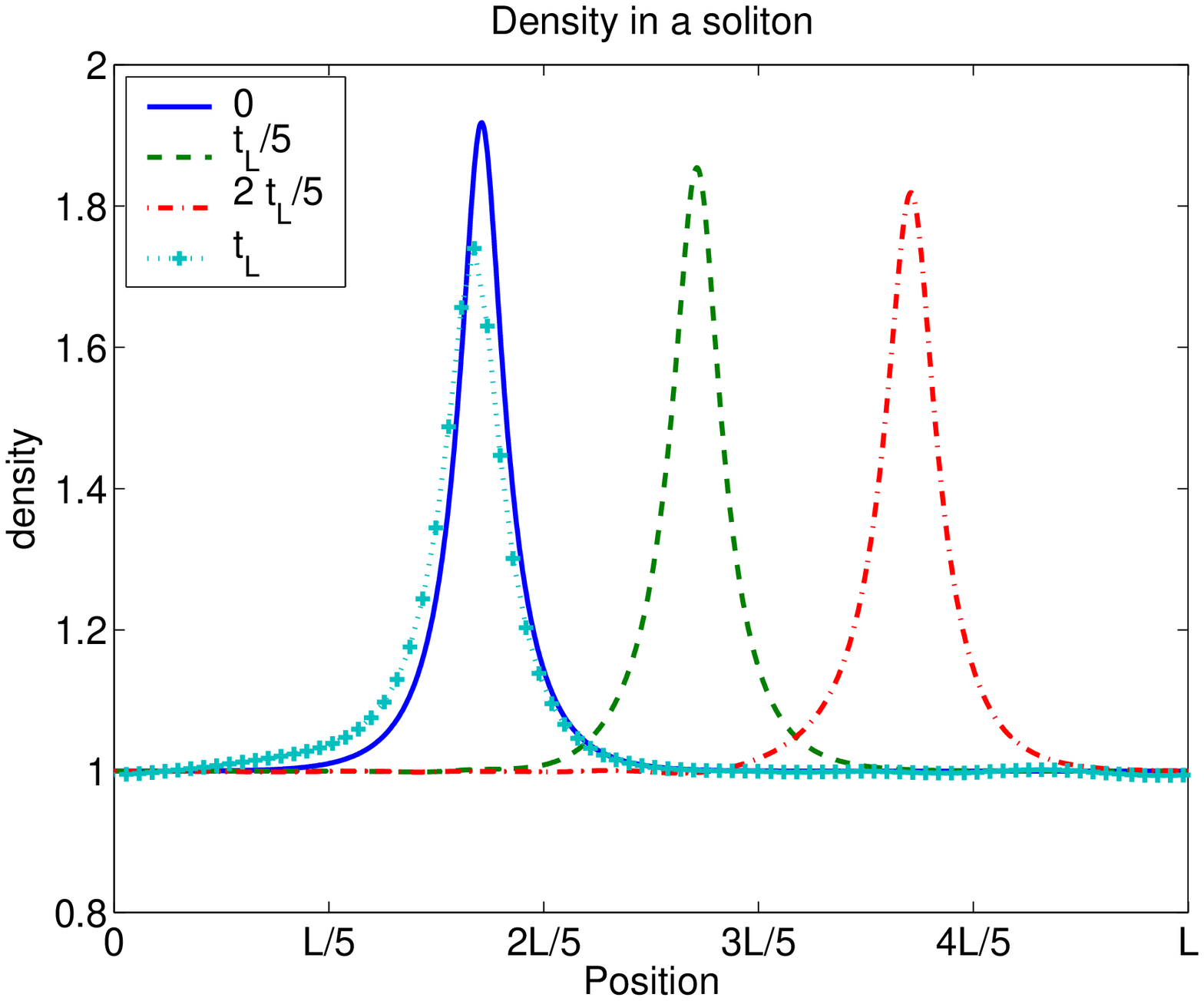}
 \end{minipage}
 \begin{minipage}[c]{.46\linewidth}
 \includegraphics[scale=0.4]{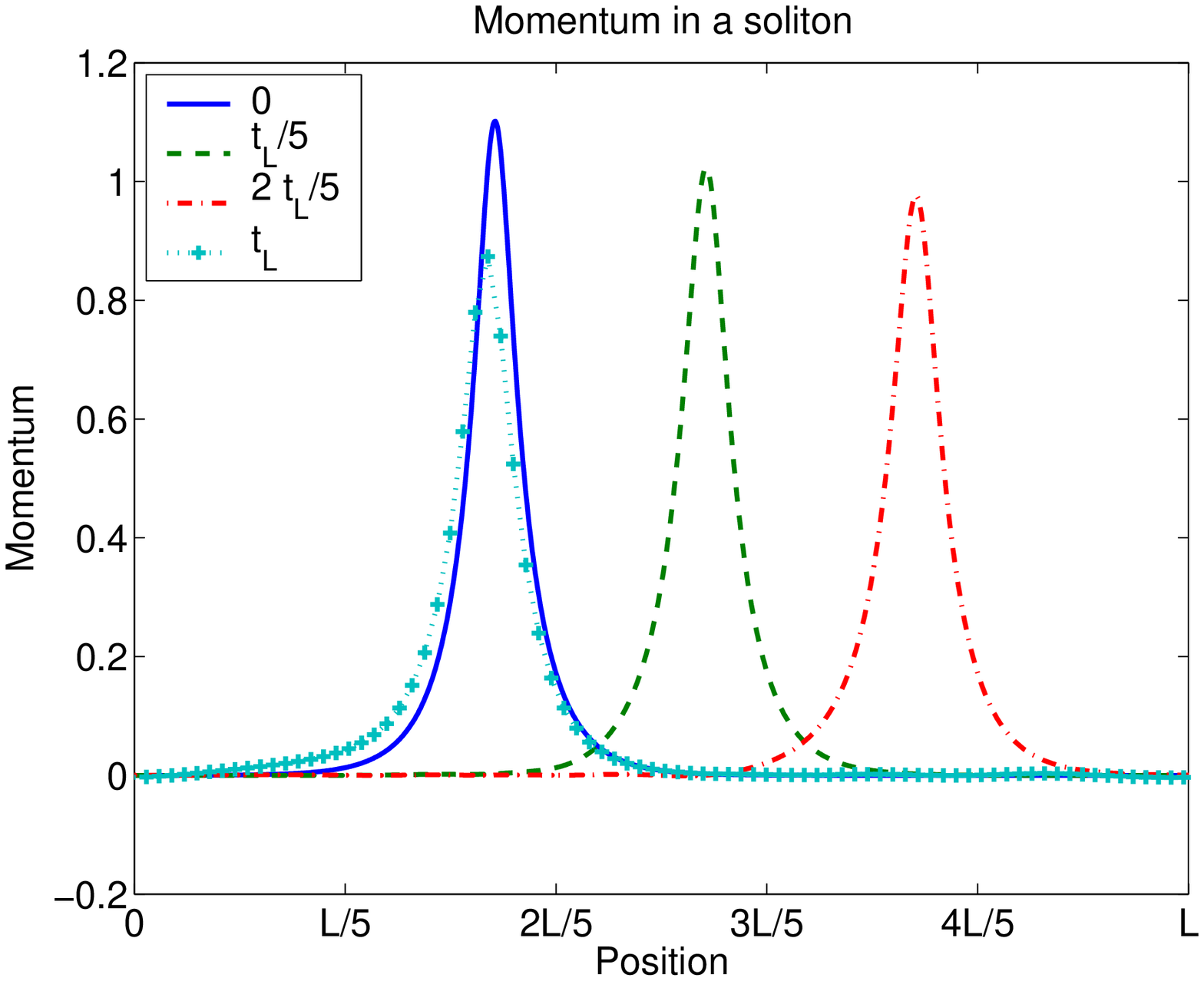}
 \end{minipage}
 \caption{\label{soliton1} Density (left) and momentum (right) as a function of space, for a soliton moving to the right,
 at time $ 0 $, $ t_{L}/5 $, $ 2 t_{L} / 5 $ and $ t_{L} $,
 computed with the REPB scheme.
 At time $ t_{L} $ the wave has travelled through the whole computational domain. 
 The damping of the density and velocity peaks are clearly visible due to the large space step, $ \Delta x = h/4 $.
 }
\end{figure}

\begin{figure}[hbtp]
 \begin{minipage}[c]{.46\linewidth}
 \includegraphics[scale=0.4]{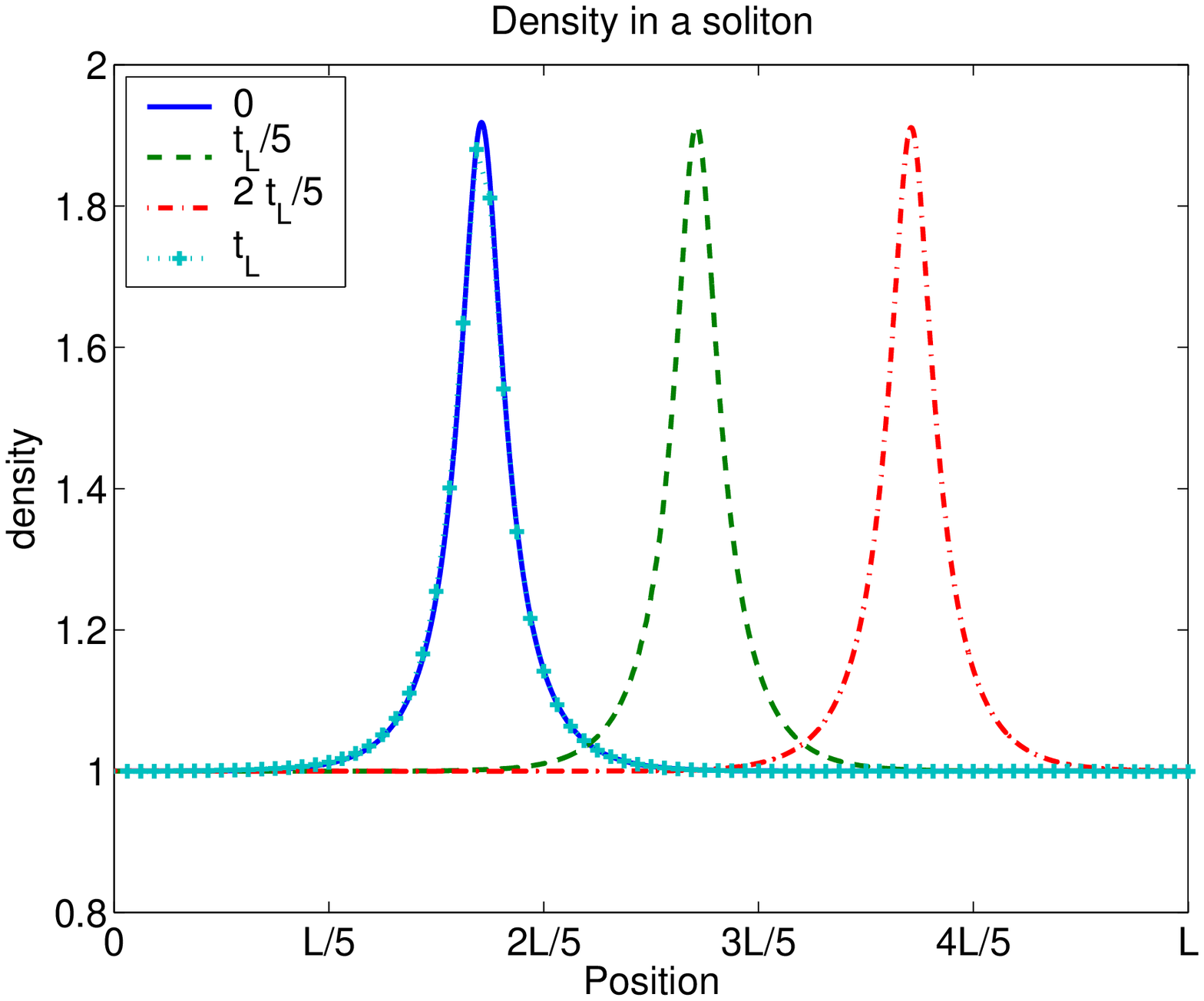}
 \end{minipage}
 \begin{minipage}[c]{.46\linewidth}
 \includegraphics[scale=0.4]{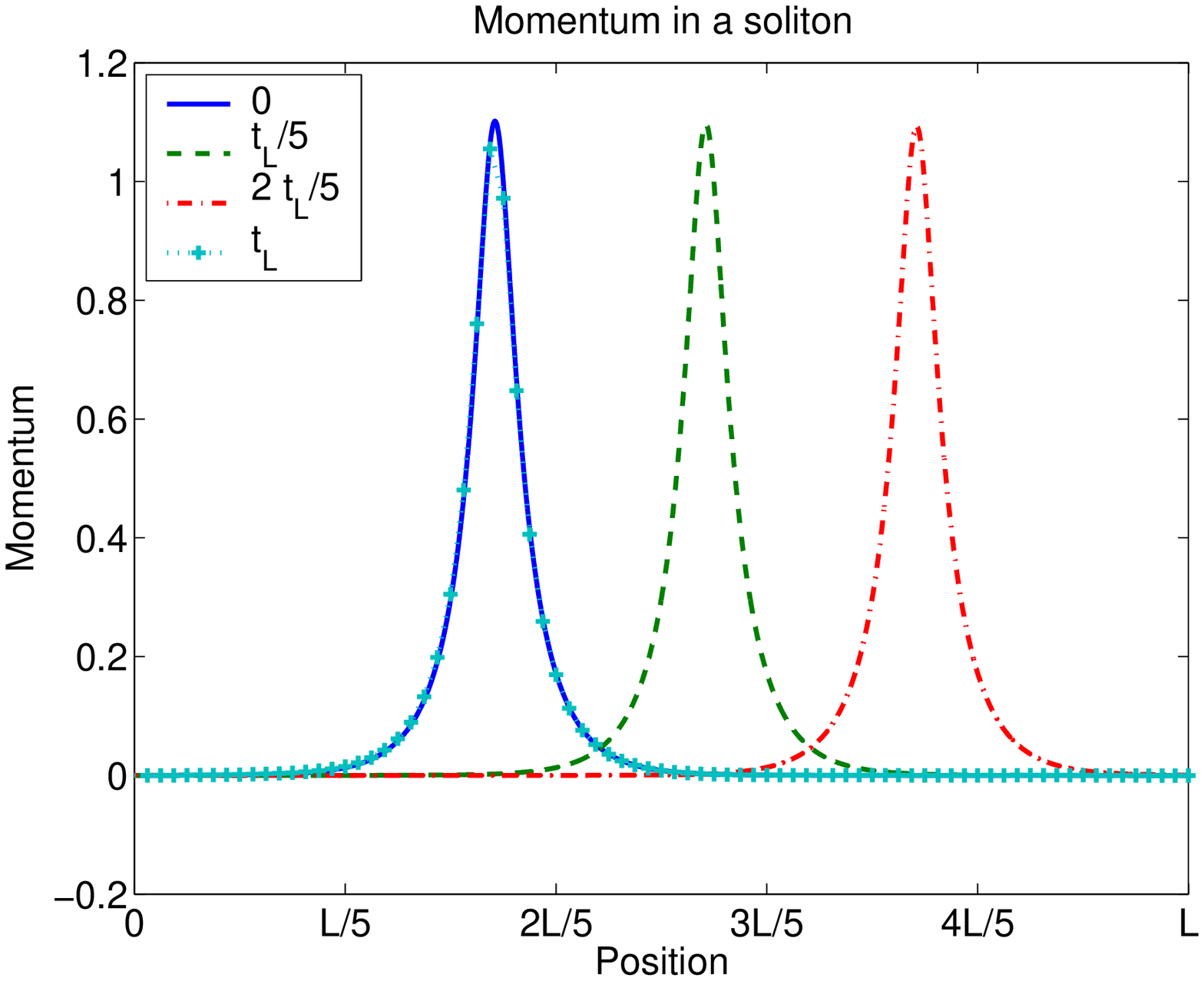}
 \end{minipage}
 \caption{\label{soliton2} Density (left) and momentum (right) as a function of space, for a soliton moving to the right,
 at time $ 0 $, $ t_{L}/5 $, $ 2 t_{L} / 5 $ and $ t_{L} $,
 computed with the REPB scheme.
 At time $ t_{L} $ the wave has travelled through the whole computational domain.
 The solution at time $ t_{L} $ can be superimposed to the solution at time $ 0 $,
 since the grid is fine enough ($\Delta x = h/64 $) and the numerical diffusion is very small.
 }
\end{figure}

The results of this comparisons are now commented. 
Figure \ref{soliton1} and \ref{soliton2} show the density and velocity in a soliton at times
 $ 0 $, $ t_L/5 $, $ 2 t_{L}/5$ and $ t_{L} $, computed with the reformulated scheme.
The shape of the initial condition is conserved with time, but the peak
 is damped and does not return to its original location after $ t_{L} $.
This effect is due to the numerical diffusion, which is inherent to the numerical method,
 and can be observed with the classical scheme as well.
In order to accurately compare the reformulated and classical schemes, one needs to perform
 a more precise study of this damping.
The forthcoming convergence study uses six grids, with a range of space steps
 from $ h $ to $ h/64 $, where $ h $ corresponds to $ 250 $ cells.
It confirms that both schemes are first order in space,
 and even if the REPB scheme suffers from a larger numerical diffusion than
 the classical one, both provide satisfactory results for this test case.

\begin{figure}[hbtp]
 \begin{minipage}[c]{.46\linewidth}
 \psfrag{varepsilon}{\footnotesize{$\varepsilon(n)$}}
 \psfrag{Space step}{\footnotesize{$\Delta x$}}
 \includegraphics[scale=0.4]{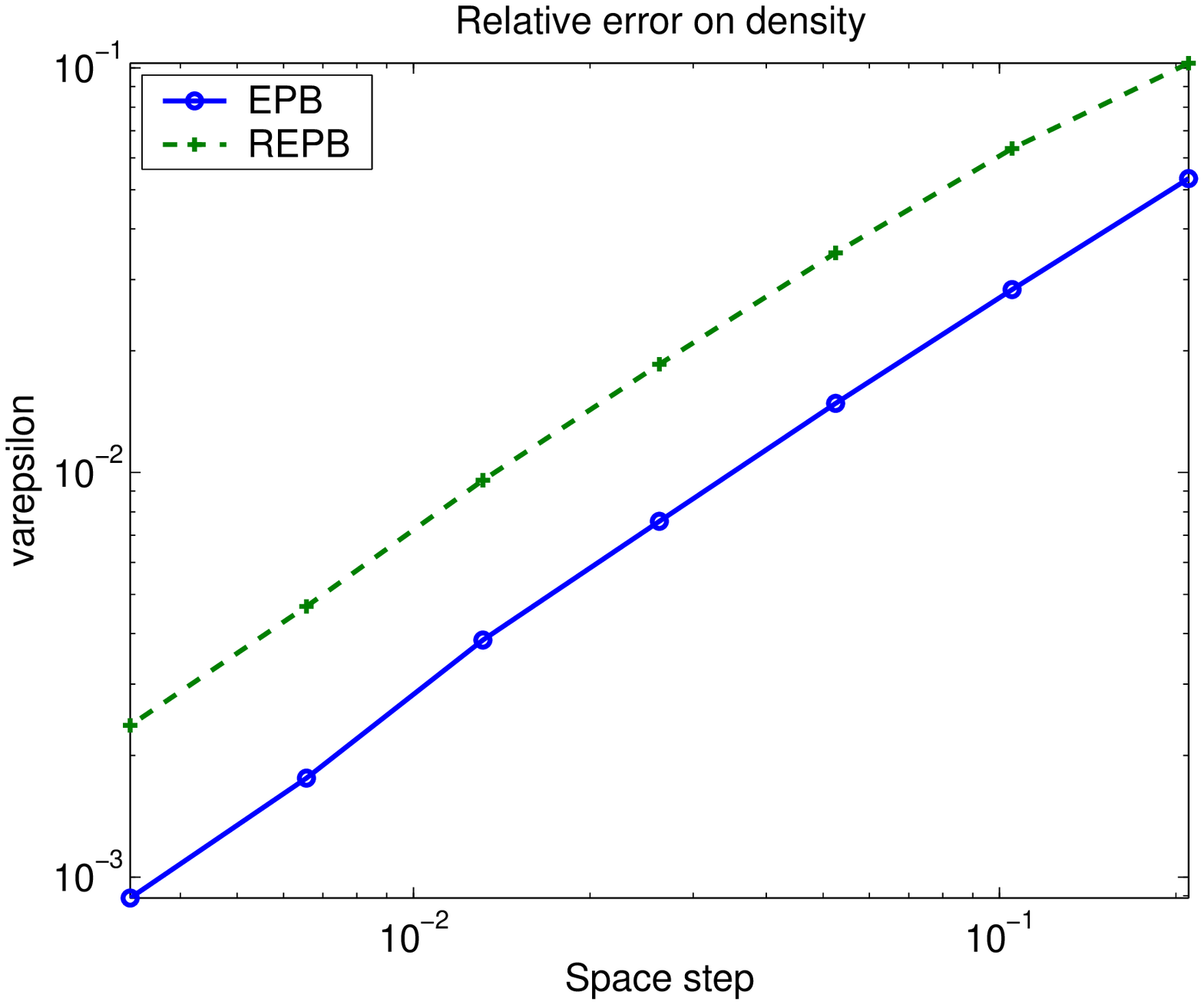}
 \end{minipage}
 \begin{minipage}[c]{.46\linewidth}
 \psfrag{varepsilon}{\footnotesize{$\varepsilon(nu)$}}
 \psfrag{Space step}{\footnotesize{$\Delta x$}}
 \includegraphics[scale=0.4]{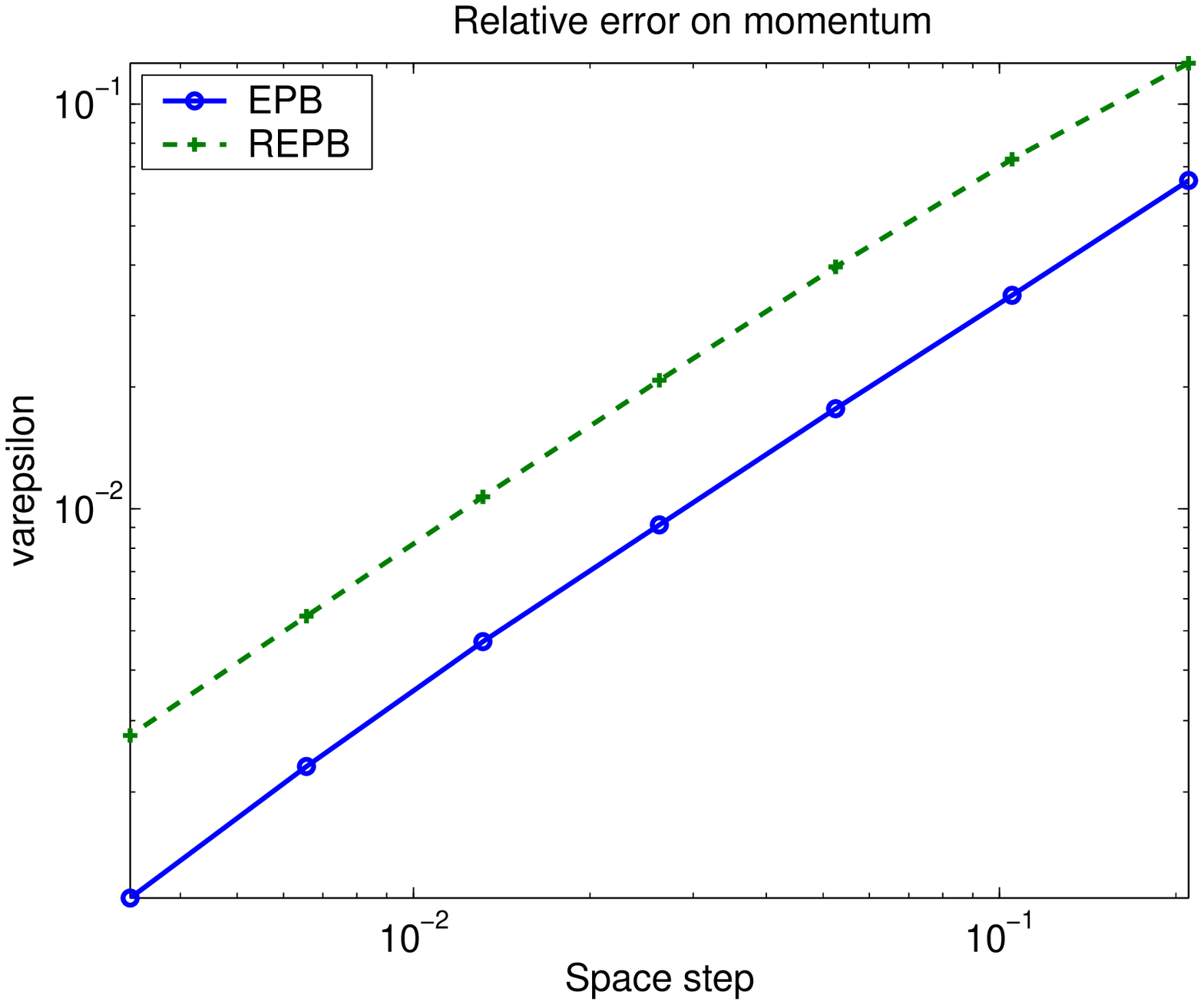}
 \end{minipage}\\
 \begin{minipage}[c]{.46\linewidth}
 \psfrag{varepsilon}{\footnotesize{$\varepsilon(\phi)$}}
 \psfrag{Space step}{\footnotesize{$\Delta x$}}
 \includegraphics[scale=0.4]{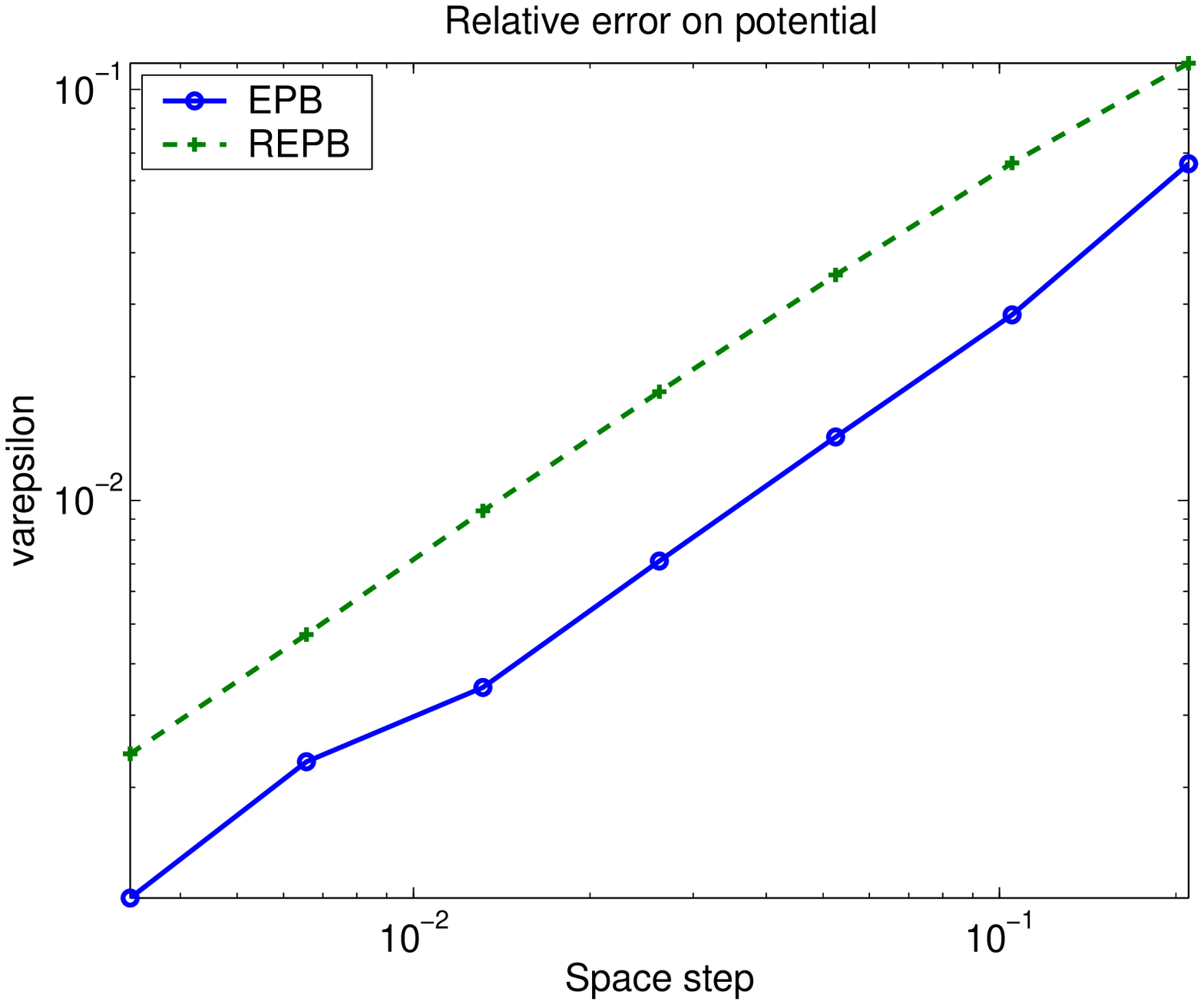}
 \end{minipage}
 \begin{minipage}[c]{.46\linewidth}
 \caption{\label{error_fig} Comparison of the error in $ L^{\infty} $ norm
 for the density, momentum and potential with various spacial mesh sizes $ \Delta x $
 for the EPB and REPB schemes.
 Both schemes are first order in space. }
 \end{minipage}
\end{figure}

The numerical convergence of the schemes in space is investigated by comparison with
 the reference solution;
at time $ t/5 $ the $ L^{\infty} $ relative error with the reference solution is computed.
For instance the density error is
$$ \varepsilon_{EPB}(n) = \frac{ \max || n_{\text{num}} - n_{\text{ref}} ||_{\infty} }{ || n_{\text{ref}} ||_{\infty} },
$$
where $ n_{\text{ref}} $ is the density of the reference solution, i.e. the initial density
 translated by $ L/5 $, and $ n_{\text{num}} $ is the density computed
 with the classical scheme.
Such errors are defined for the reformulated scheme and for the momentum
 and potential.
The schemes are tested on six grids, made of $ 250,500,1000,2000,4000, 8000 $ and $ 16000 $ cells.
Table \ref{error_tab} and figure \ref{error_fig} confirm that
 both numerical schemes are first order in space.
The REPB scheme suffers from a larger numerical dissipation than the EPB scheme.

\begin{table}[hbtp]
\begin{tabular}{|c|c|c||c|c||c|c|}
\hline
N & $ \varepsilon_{EPB}(n) $ & $ \varepsilon_{REPB}(n) $
  & $ \varepsilon_{EPB}(nu) $ & $ \varepsilon_{REPB}(nu) $
  & $ \varepsilon_{EPB}(\phi) $ & $ \varepsilon_{REPB}(\phi) $ \\
\hline
\hline
$h$    & $ 0.053 $&$ 0.102 $&$ 0.111 $&$ 0.215 $&$ 0.066 $&$ 0.116 $ \\
$h/2$  & $ 0.028 $&$ 0.060 $&$ 0.060 $&$ 0.131 $&$ 0.028 $&$ 0.066 $ \\
$h/4$  & $  0.014 $&$ 0.035 $&$ 0.032 $&$ 0.073 $&$ 0.014 $&$ 0.035 $ \\
$h/8$  & $  7.6\times 10^{-3} $&$ 18.5\times 10^{-3} $&$ 1.64\times 10^{-2} $&$
         3.87\times 10^{-2} $&$ 7.12\times 10^{-3} $&$ 18.4\times 10^{-3} $ \\
$h/16$ & $ 3.86\times 10^{-3} $&$ 9.57\times 10^{-3} $&$ 8.40\times 10^{-3} $&$
         2.00\times 10^{-2} $&$ 3.50\times 10^{-3} $&$ 9.42\times 10^{-3} $ \\
$h/32$ & $ 1.76\times 10^{-3} $&$ 4.67\times 10^{-3} $&$ 3.83\times 10^{-3} $&$
         1.00\times 10^{-2} $&$ 2.31\times 10^{-3} $&$ 4.71\times 10^{-3} $ \\
$h/64$ & $ 8.89\times 10^{-4} $&$ 2.37\times 10^{-3} $&$ 1.97\times 10^{-3} $&$
         5.09\times 10^{-3} $&$ 1.07\times 10^{-3} $&$ 2.41\times 10^{-3} $ \\
\hline
\end{tabular}
\caption{\label{error_tab} Comparison of the error in $ L^{\infty} $ norm
 for the density, momentum and potential with various space grid sizes
 using the EPB and REPB schemes.
The computation of the error is made with the analytical solution for
 the soliton test case}
\end{table}

\begin{figure}[hbtp]
 \includegraphics[scale=0.4]{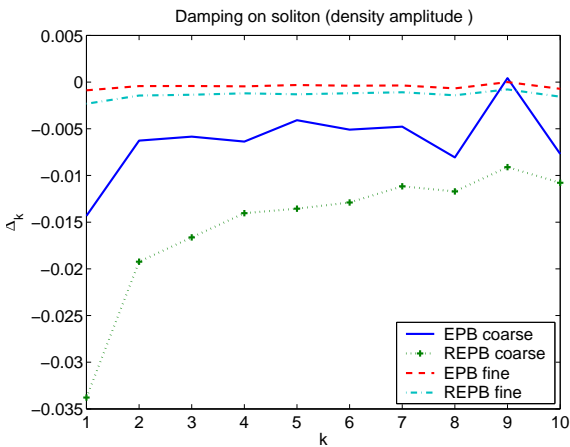}
 \caption{\label{damping_fig} 
 Comparison of the damping rates for the density amplitude of the soliton for the EPB and REPB schemes on two different spatial mesh sizes:
 a coarser grid with $ 1000 $ cells and a finer grid with $ 16000 $ cells are used for both schemes.}
\end{figure}

The numerical dissipation of the schemes can be measured in another way.
Indeed, the amplitude of the numerical soliton is damped with time.
The following study quantifies the damping rates of both the EPB and REPB schemes.
This study is performed on a long time simulation: its final time is $ 2 t_{L} $.
This study compares the damping rate of the density amplitude over one time increment $ \Delta t $
 (not related with the numerical time step), at different times of the simulation.
The time increment $ \Delta t $ is $ t_{L}/5 $.
The density amplitude $ n_{\max} ((k+1)\Delta t)$ is compared to the density amplitude at the previous time increment $ n_{\max} (k \Delta t)$
 for $ 0 \leq k \leq 9 $.
We measure the decrement of the amplitude $ \Delta_{k} $ between $ k\Delta t $ and $ (k+1)\Delta t $ as follows:
\begin{equation}
\Delta_{k} = \frac{1}{\Delta t} 
\left| \ln \left( \frac{ n_{\max}((k+1)\Delta t)}{n_{\max}(k\Delta t)} \right) \right|
\end{equation}
The damping rate of the wave amplitude appears on figure \ref{damping_fig} for both the EPB and REPB schemes.
Two spatial grids are used for this comparison: a coarse grid with $ 1000 $ cells and a fine grid with $ 16000 $ cells.
The REPB scheme shows a larger dissipation than the classical one.
The evolution of the damping rates with simulation time is similar for the two schemes :
on a coarse grid the peak is damped faster at the beginning of the simulation.

\subsection{Riemann problem}
\label{subsec_riemann}

The previous test case, in which the domain size and Debye length are of the same order of magnitude, was designed for the dispersive regime $\lambda = O(1)$. The present test case aims at investigating how the schemes perform in the hydrodynamic (or quasi-neutral) regime $\lambda \ll 1$. Therefore, in this test case, values of $\lambda$ ranging from small to very small small are used.
When $\lambda$ is small, the REPB model explicitly appears as a perturbation of the ICE model, which is a strongly hyperbolic and conservative model. Therefore the REPB-based scheme should be accurate in the hydrodynamic regime, in particular for the computation of solutions involving discontinuities.
By contrast, the hydrodynamic part of the EPB model is a pressureless gas dynamics model, which is not strictly hyperbolic, and is thus weakly unstable.
For this reason, the EPB-based scheme is expected to be less accurate in the small Debye length limit. The present test problem aims at testing the validity of these predictions. 

The test case is a shock tube problem involving two outgoing shock waves.
The initial density is a constant, while the initial velocity has a jump at $ x = 0 $ between the constants $ u_{L} = +1 $ and $ u_{R} = -1$.
The density in the intermediate state of the Riemann problem is larger and the velocity is zero. Two outgoing shock waves appear on each side of this intermediate state.
The computational domain is $ [-0.2;0.2] $. The dimensionless Debye
 length $ \lambda $ varies from $ 10^{-2} $ to $ 10^{-4} $.

The value $ \lambda = 10^{-2} $ is large enough for singularities near the shock waves to appear, due to the dispersive nature of the  Euler-Poisson-Boltzmann model. In the semi-classical setting, the framework of multi-valued solutions \cite{Liu_Wang_1,Liu_Wang_2} can be used to explain the qualitative features of the classical solutions. Indeed, classical solutions keep a signature of these underlying multi-valued solutions, in the form of singularities (when the number of branches changes) and oscillations (when several branches co-exist and the solution 'oscillates' between these branches). In the $\lambda \to 0$ limit, the entropic solutions of the ICE model capture the average value of these oscillations, but not the details of them. As we will see, the EPB-based scheme keeps better track of these oscillations but when the mesh size does not resolve the Debye length, the oscillations become mesh-dependent. On the other hand, the REPB-based scheme directly provides the entropic solution of the limiting ICE model and is better suited to capture the average value of these oscillations (i.e. the weak limit of the solutions of the EPB model when $\lambda \to 0$).

\begin{figure}[hbtp]
 \begin{minipage}[c]{.46\linewidth}
 \includegraphics[scale=0.4]{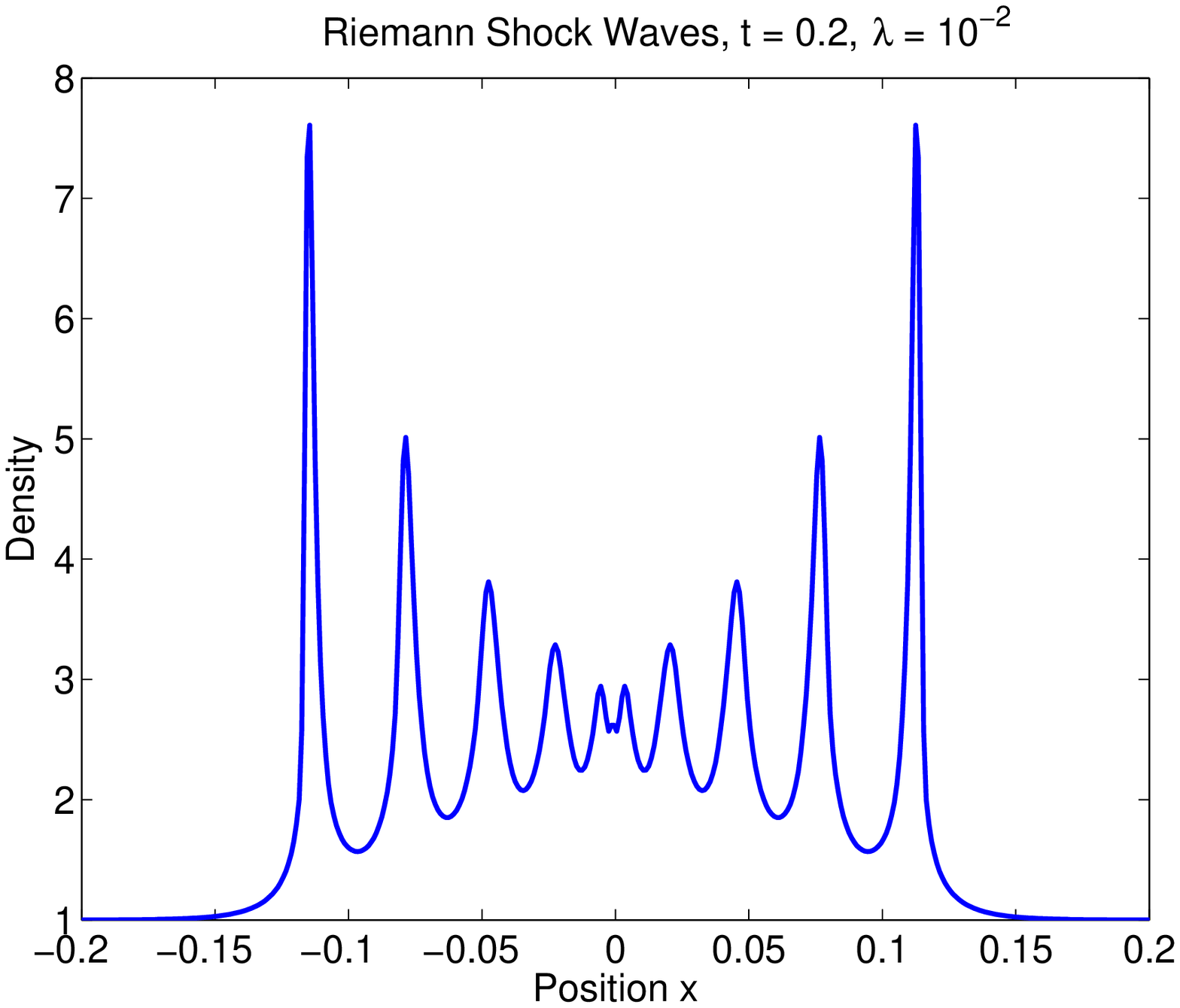}
 \end{minipage}
 \begin{minipage}[c]{.46\linewidth}
 \psfrag{varepsilon}{\footnotesize{$\varepsilon(nu)$}}
 \psfrag{Space step}{\footnotesize{$\Delta x$}}
 \includegraphics[scale=0.4]{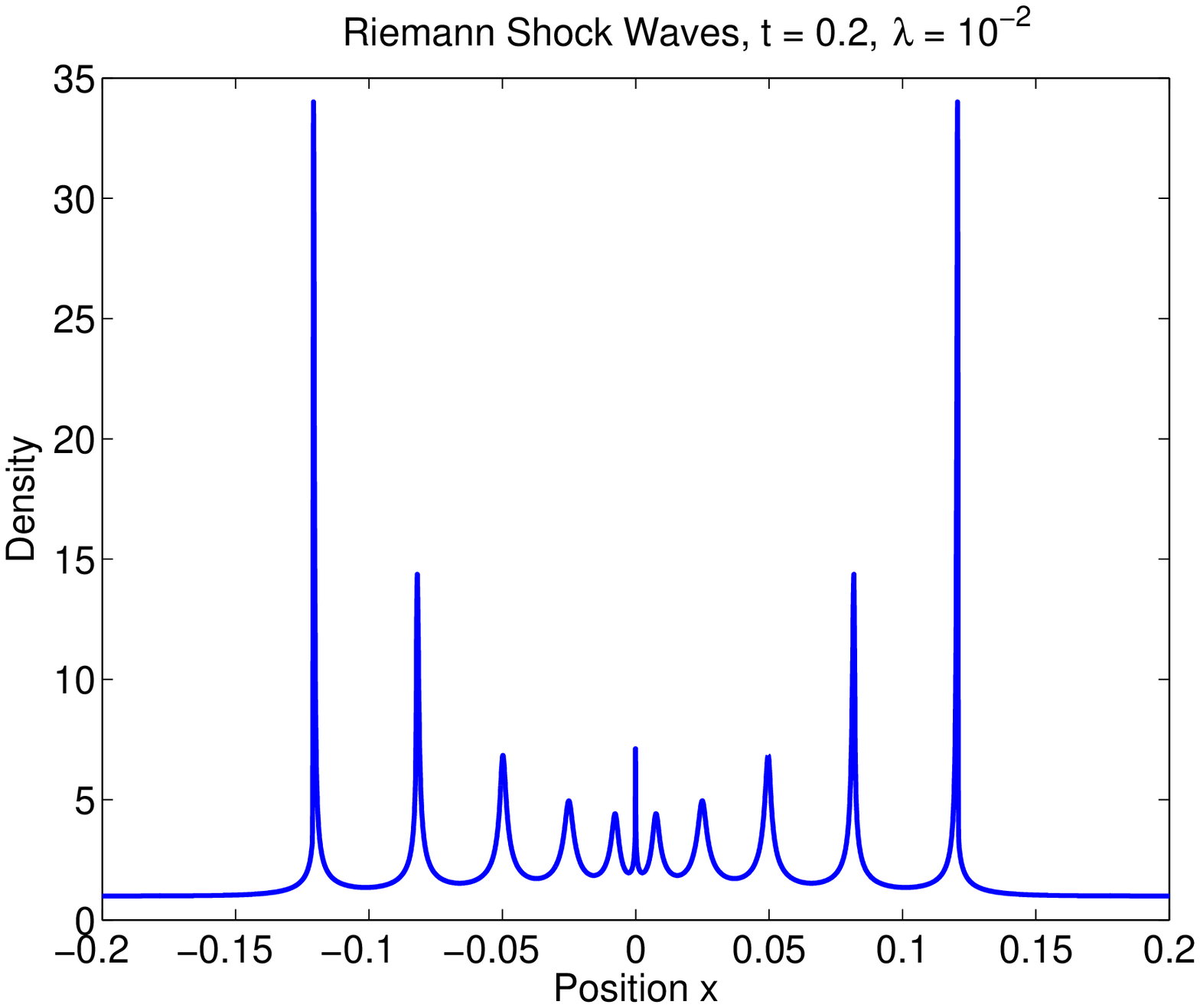}
 \end{minipage}
\caption{\label{densities}Density as a function of space for the Riemann shock-wave problem
 at time $ t = 0.2 $ with dimensionless Debye length $ \lambda = 10^{-2} $.
Left: classical scheme with a coarse grid ($ N = 2000 $). Right: 
classical scheme with a fine grid ($ N = 32000 $).
Density peaks are sharpened when the mesh size decreases}
 \begin{minipage}[c]{.46\linewidth}
 \psfrag{varepsilon}{\footnotesize{$\varepsilon(\phi)$}}
 \psfrag{Space step}{\footnotesize{$\Delta x$}}
 \includegraphics[scale=0.4]{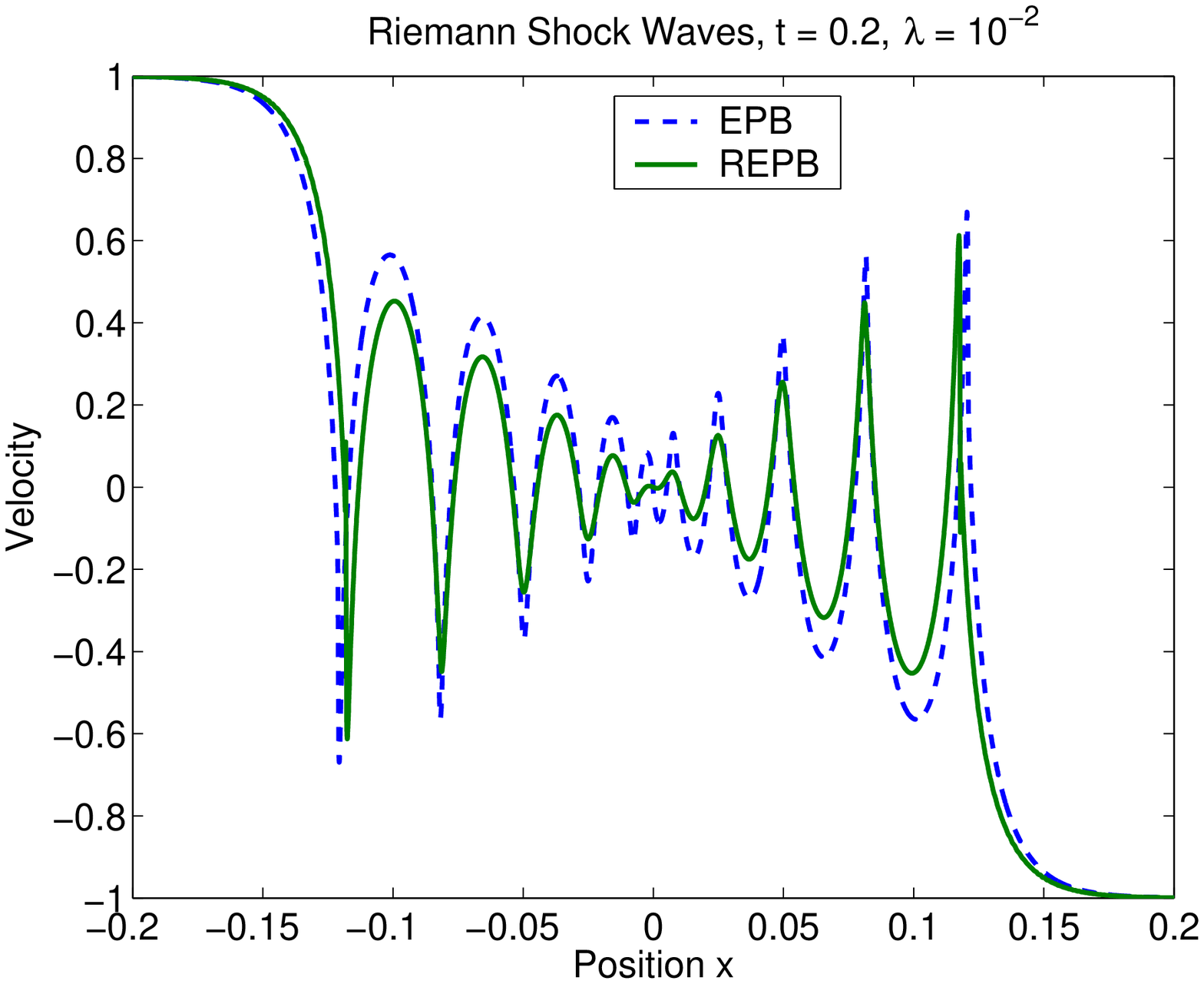}
 \end{minipage}
 \begin{minipage}[c]{.46\linewidth}
\includegraphics[scale=0.4]{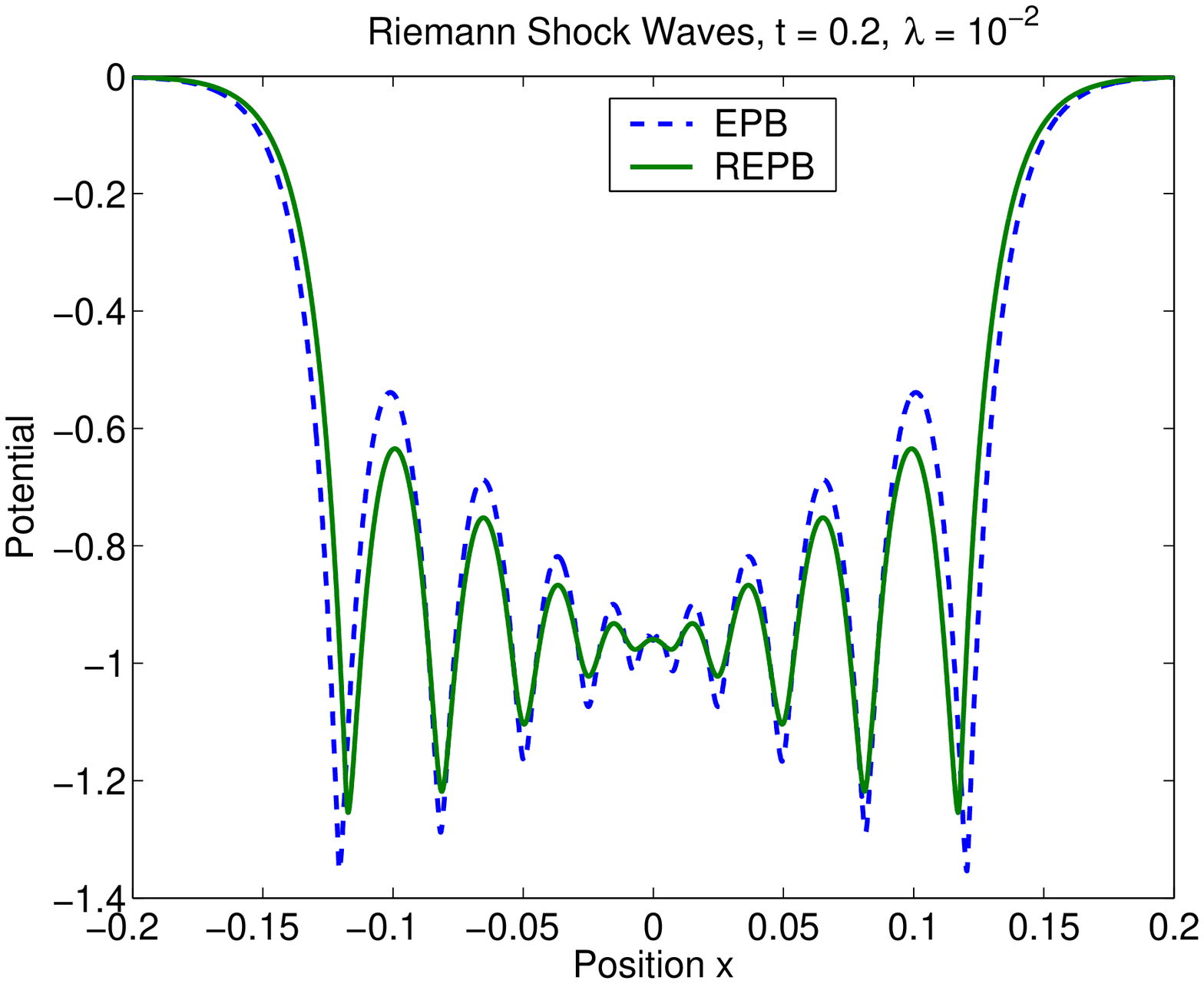}
 \end{minipage}
\caption{\label{vit_pot}Velocity (left) and potential (right) as a function of space for the Riemann shock-wave problem
 at time $ t = 0.2 $ with dimensionless Debye length $ \lambda = 10^{-2} $.
EPB-based (dashed line) and REPB-based (solid line) schemes on a fine grid ($ N = 32000 $). }
\end{figure}

Fig. \ref{densities} confirms that the numerical solution provided by the EPB-based scheme in under-resolved situation is mesh-dependent. It displays the density computed by the EPB-based schemes for two different mesh sizes. The density peaks computed on the finer mesh are much finer and higher than those computed on the coarse mesh. By contrast, the velocity and potential remain finite regardless of the mesh size (not displayed).
Fig. \ref{vit_pot} shows the velocity and potential computed with the EPB-based scheme compared to the REPB-based scheme.
Like in the soliton test case, the REPB-based scheme shows a slightly larger numerical dissipation.

\begin{figure}[hbtp]
 \begin{minipage}[c]{.46\linewidth}
 \includegraphics[scale=0.38]{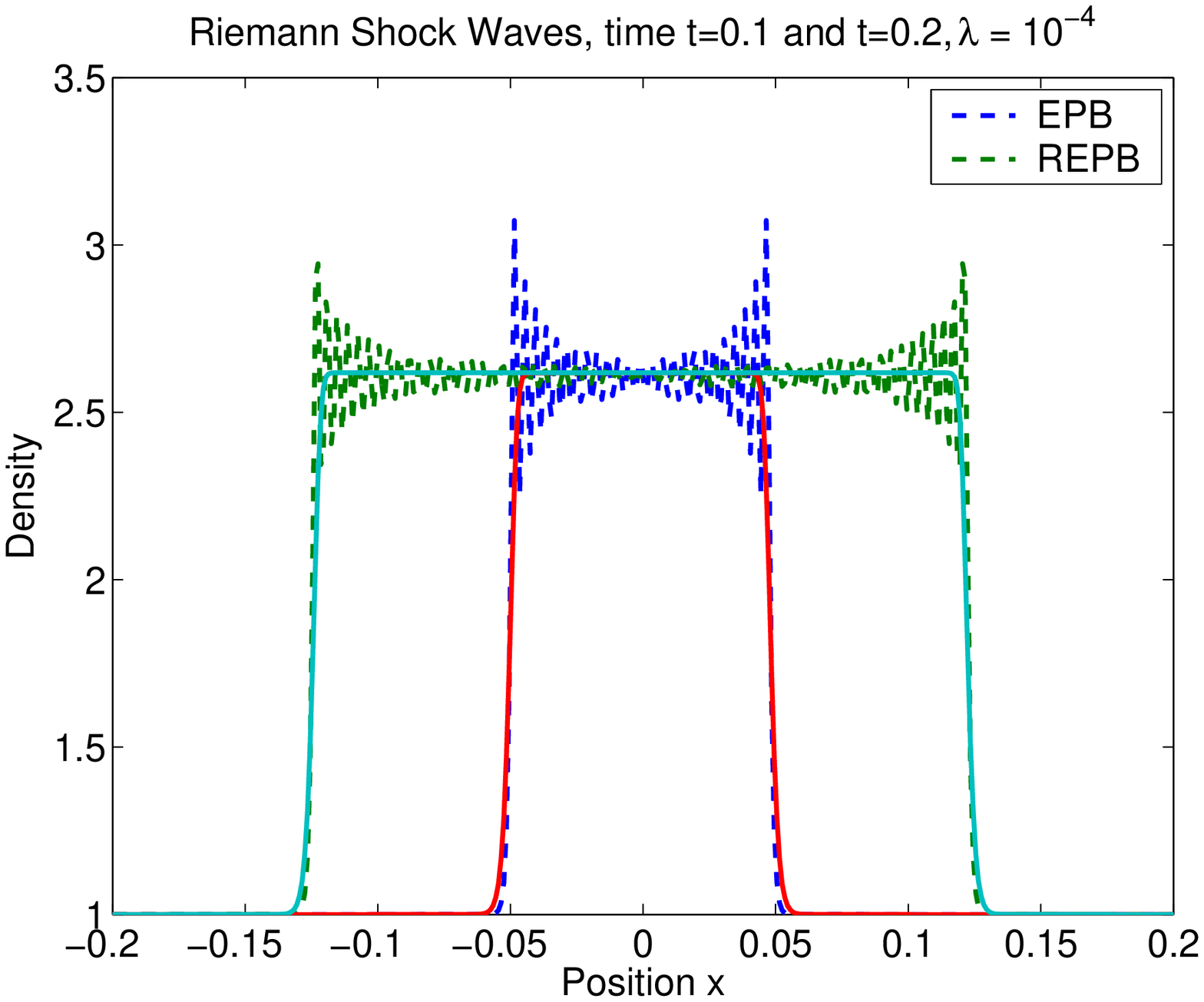}
 \end{minipage}
 \begin{minipage}[c]{.46\linewidth}
 \psfrag{varepsilon}{\footnotesize{$\varepsilon(nu)$}}
 \psfrag{Space step}{\footnotesize{$\Delta x$}}
 \includegraphics[scale=0.38]{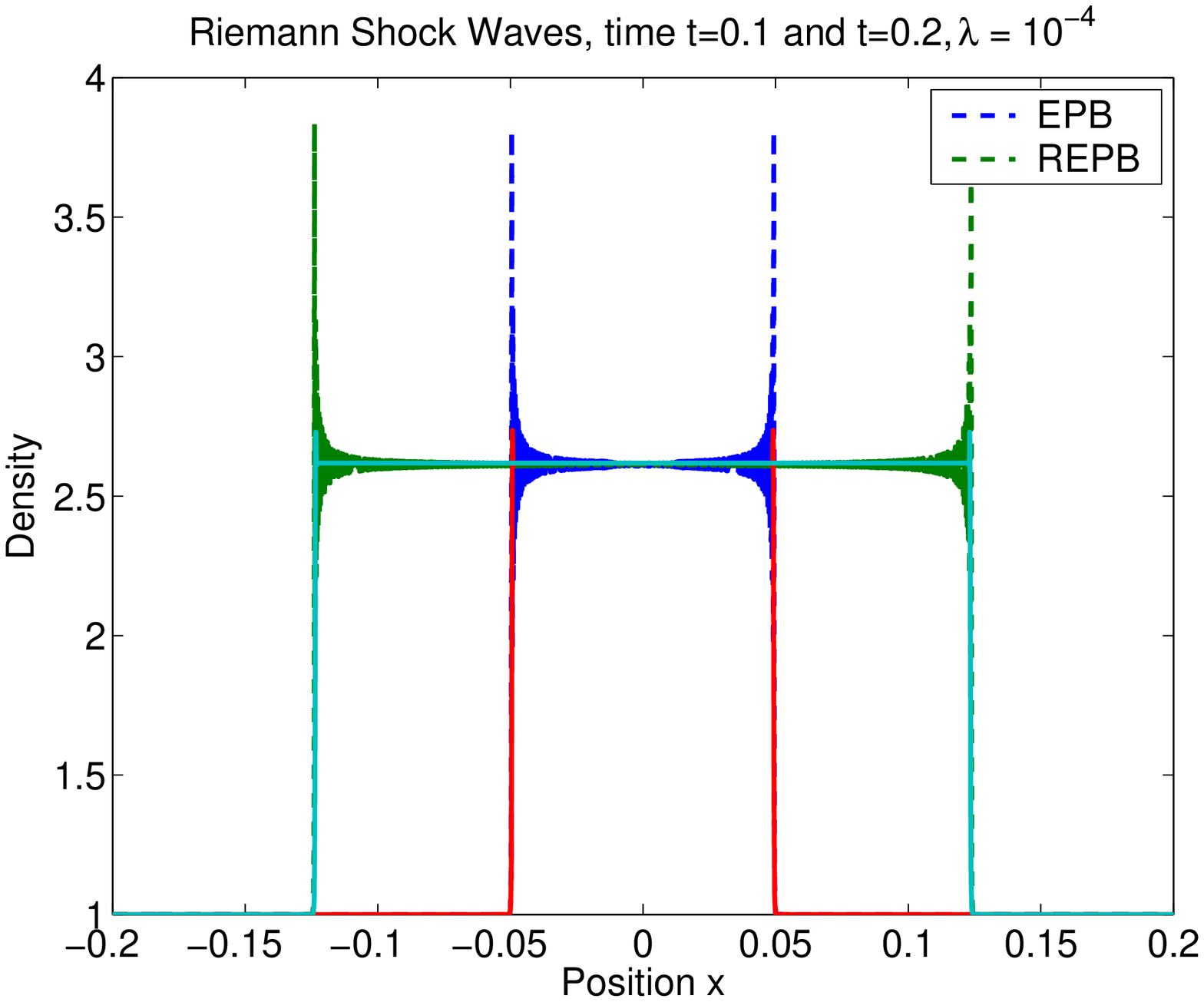}
 \end{minipage}
\caption{\label{4_densities}Density as a function of space for the Riemann shock-wave problem
 at time $ t = 0.1 $ and $ t = 0.2 $ with dimensionless Debye length $ \lambda = 10^{-4} $. Left: 
EPB and REPB-based schemes on a coarse grid ($ N = 2000 $). Right: 
same schemes on a fine grid ($ N = 32000 $).}
\end{figure}

The value $ \lambda = 10^{-4} $ is small enough to observe the hydrodynamic regime.
Both schemes show an accurate determination of the shock speed but the EPB scheme leads to spurious oscillations in the neighborhood of the shock.
These oscillations are shown on figure \ref{4_densities}.
They are mesh dependent and occur even with smaller dimensionless Debye length $ \lambda $.
Whatever the mesh size is the solution computed with the reformulated REPB scheme does not present such oscillations. 

This test case shows that the EPB and REPB-based schemes have a very different behavior when $\lambda \ll 1$. In such under-resolved situations, while the EPB scheme presents mesh-dependent oscillations of finite amplitude, the REPB-based scheme provides an accurate approximation of the entropic solution of the limiting ICE model. As a conclusion, the REPB-based scheme should be preferred for the $\lambda \ll 1$ regime.

\subsection{Dispersive solutions test-cases}
\label{subsec_HLtest}

\subsubsection{Description}
\label{subsubsec_dispersive_test_description}

In this section, we present test-cases which are inspired from \cite{Liu_Wang_2}. In \cite{Liu_Wang_2}, the goal was to explore the computation of multi-valued solutions in the semi-classical setting by means of the level-set method. Here we consider only classical solutions. 
The first test case is referred to as a five branch test case (because it corresponds to the occurrence of a five branch multi-valued solutions in the semi-classical setting).
The second test is a seven branch test case (for the same reason).
The initial densities of both test-cases are bumps with a gaussian shape.
This shape leads to a potential well which generates an induced electric field. This electric field in turn contracts the density bump leading to a positive feedback amplification.
This effect can be further amplified by setting up appropriate initial velocities.
In the dispersive regime the amplification of the density peak can lead to singularities, whereas in the hydrodynamic regime
 the density spreads out in the entire computational domain and its profile remains smooth.

The emergence of singularities has been investigated by Liu and Wang \cite{Liu_Wang_1, Liu_Wang_2}.
The authors compute multi-valued solution for similar test cases but with the standard Poisson equation (without the exponential term coming from the Boltzmann relation) that allows the computation of analytical solutions.
Here, no analytical solution is available but the two schemes (EPB and REPB) are tested one against each other and against a numerically computed reference solution on a very fine mesh. Due to the singularities appearing when $ \lambda = 1 $ the numerical error is computed with the potential, which remains finite in every situation. The tests are run with $ \lambda = 1 $ and $ \lambda = 10^{-2} $ in order to explore both the dispersive and hydrodynamic regimes. The results would be similar if $\lambda$ was further reduced.

\subsubsection{Five-branch solution}
\label{subsubsec_test1}

\begin{figure}[hbtp]
 \begin{minipage}[c]{.46\linewidth}
 \includegraphics[scale=0.4]{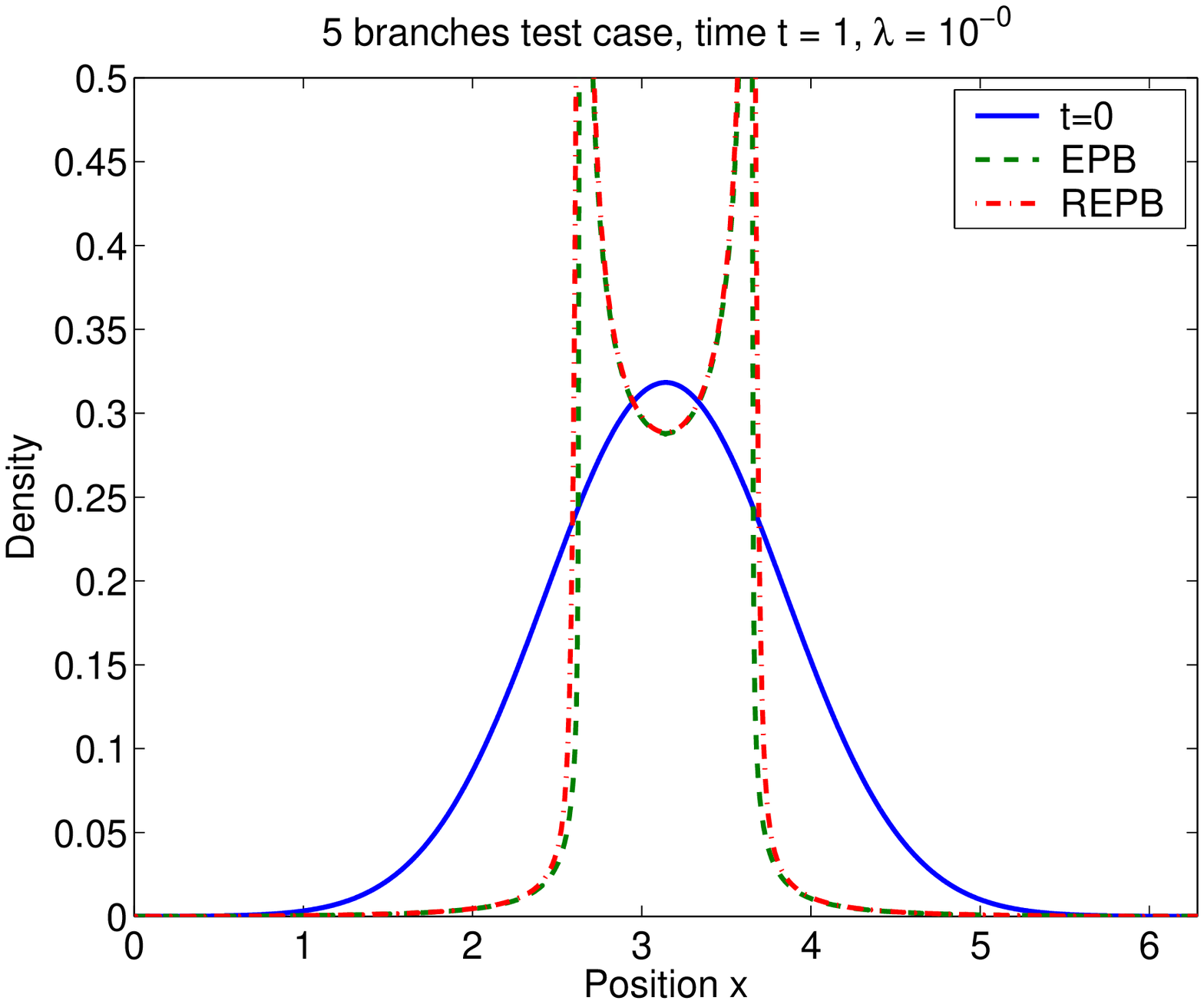}
 \end{minipage}
 \begin{minipage}[c]{.46\linewidth}
 \psfrag{varepsilon}{\footnotesize{$\varepsilon(nu)$}}
 \psfrag{Space step}{\footnotesize{$\Delta x$}}
 \includegraphics[scale=0.4]{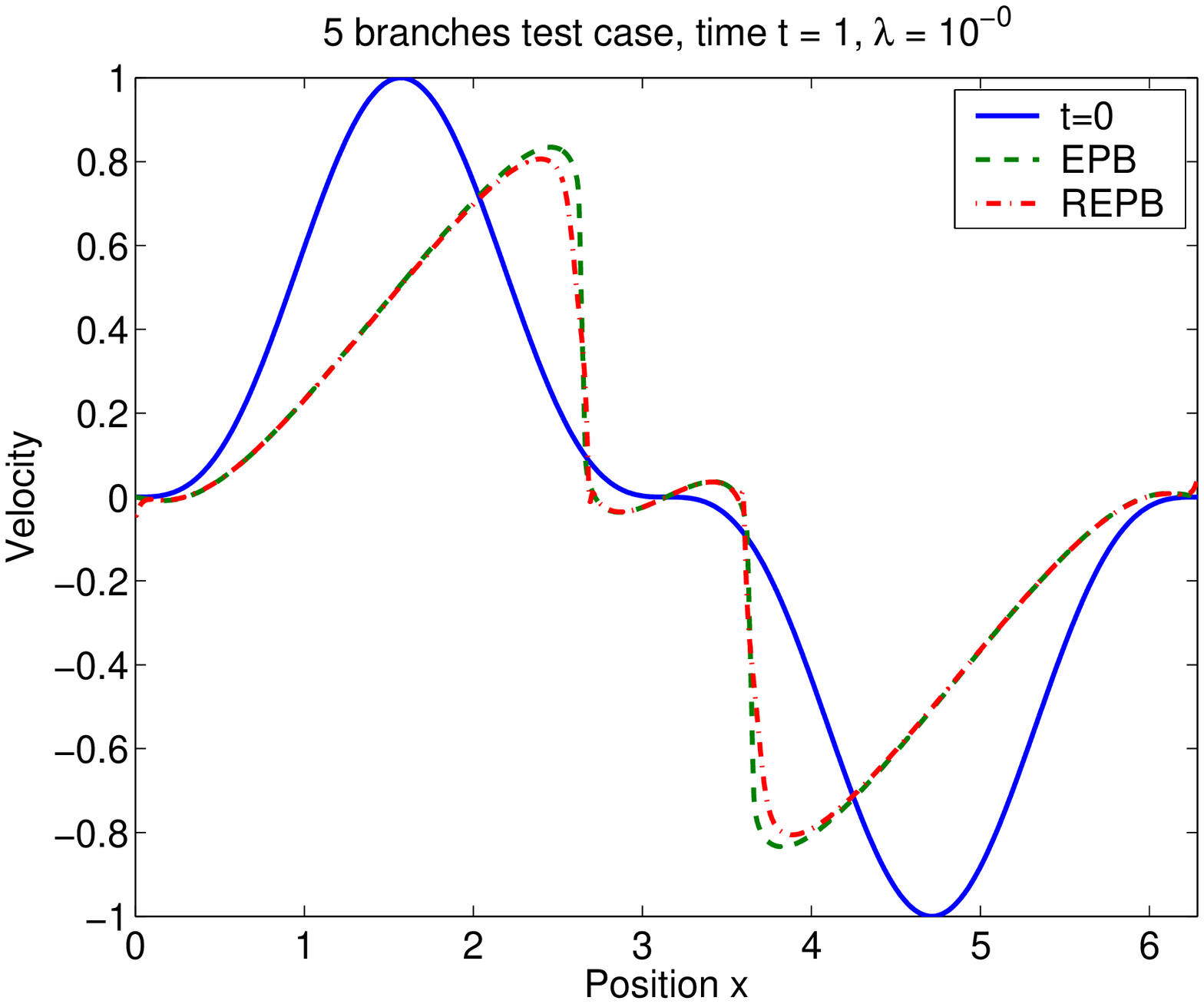}
 \end{minipage}
\caption{\label{ct5b_0_dv}Density (left) and velocity (right) as a function of space for the five-branch test case.
The dimensionless Debye length $ \lambda $ is $ 1 $. Numerical solutions at time $ t = 1 $ are computed with the EPB and REPB-based
 schemes on a grid with $ 2000 $ cells.
}
 \begin{minipage}[c]{.46\linewidth}
 \psfrag{varepsilon}{\footnotesize{$\varepsilon(\phi)$}}
 \psfrag{Space step}{\footnotesize{$\Delta x$}}
 \includegraphics[scale=0.4]{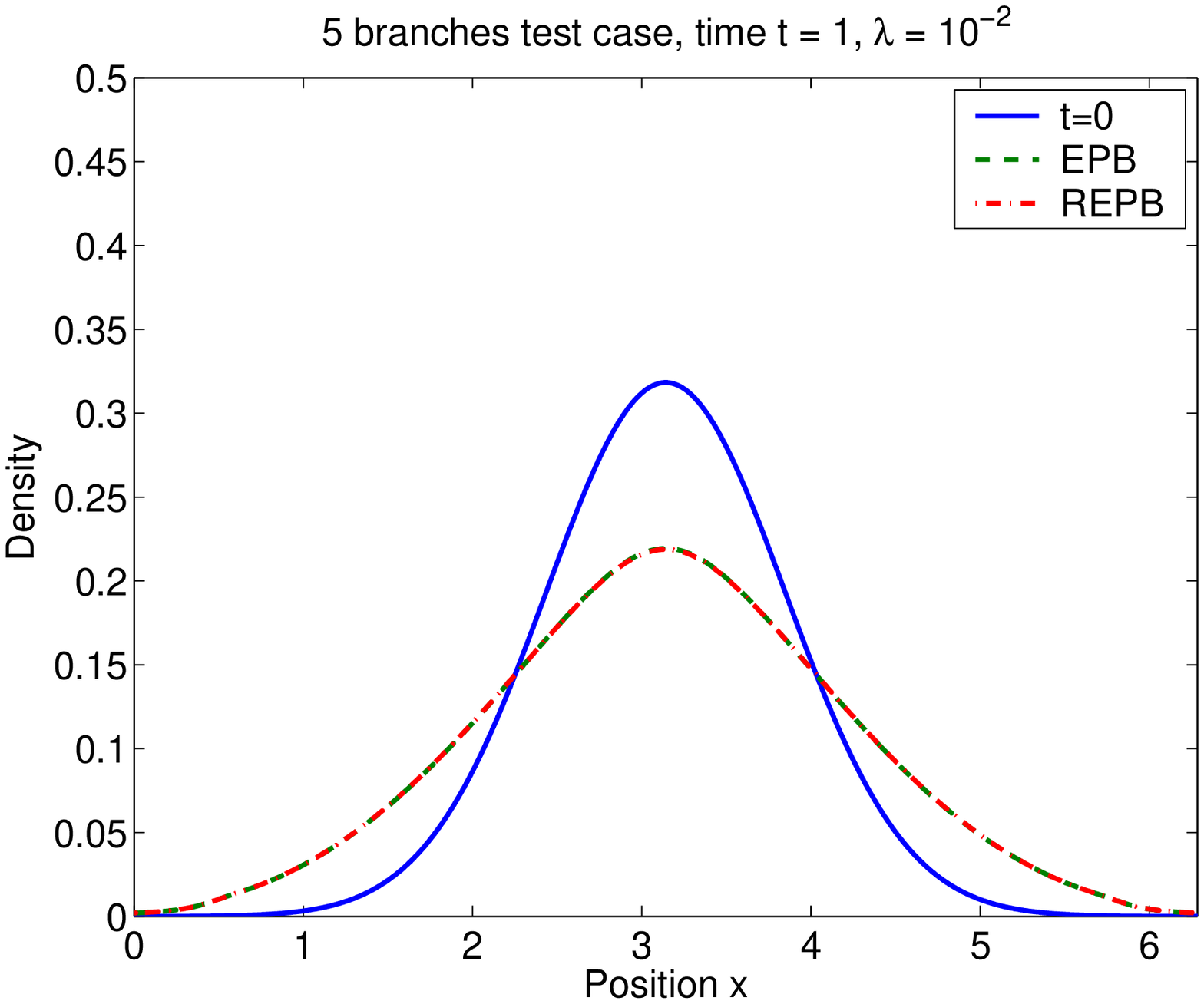}
 \end{minipage}
 \begin{minipage}[c]{.46\linewidth}
\includegraphics[scale=0.4]{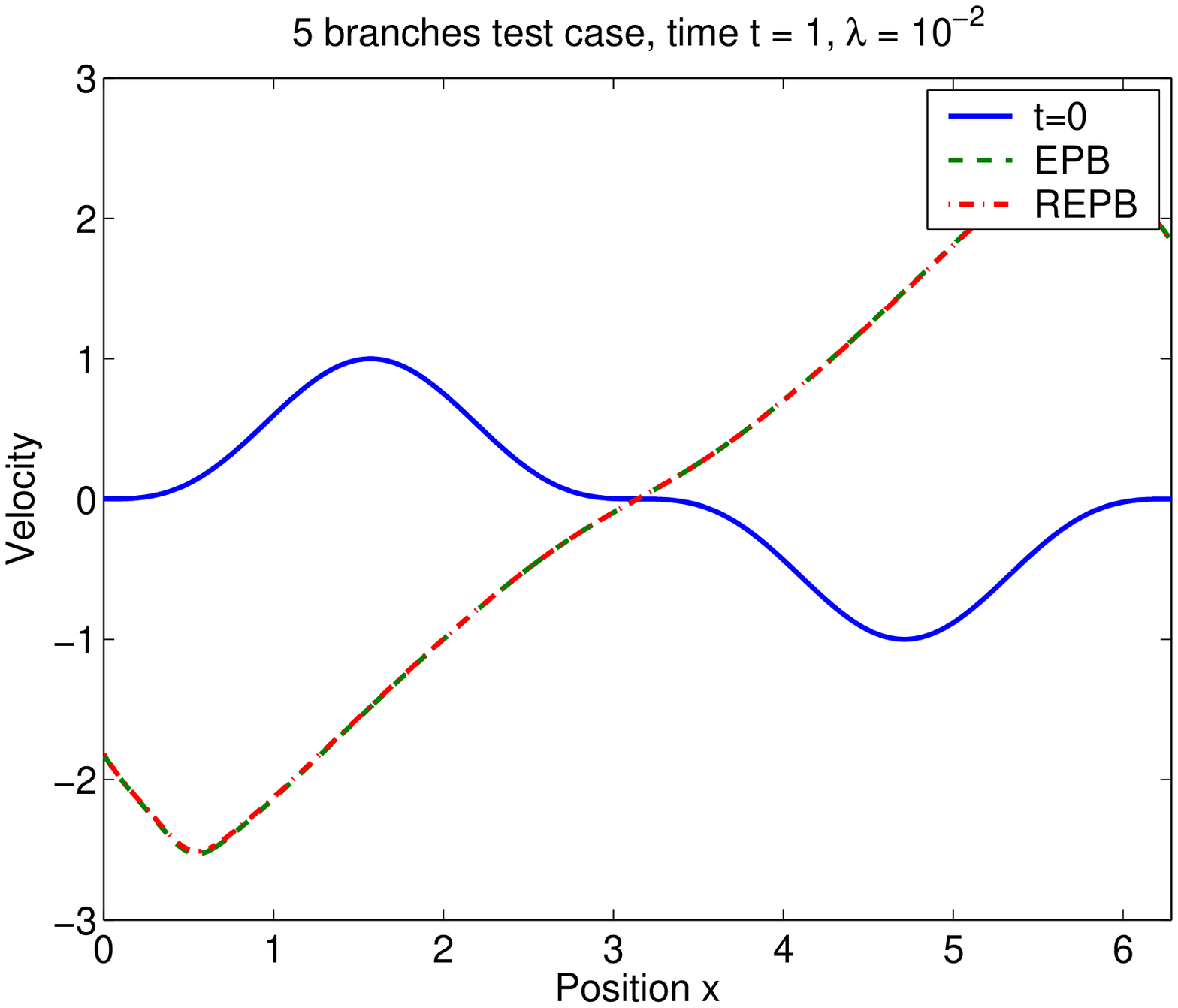}
 \end{minipage}
\caption{\label{ct5b_2_dv}Density (left) and velocity (right) as a function of space for the five-branch test case.
The dimensionless Debye length $ \lambda $ is $ 10^{-2} $. Numerical solutions at time $ t = 1 $ are computed with the EPB and REPB-based
 schemes on a grid with $ 2000 $ cells.}
\end{figure}

The initial condition for this first test-case is given by
\begin{eqnarray*}
& & n_{0}  = \frac{1}{\pi} e^{ -(x-\pi)^2 }, \\
& & u_{0}  = \sin^{3}x.
\end{eqnarray*}
The computational domain is $ [0,2 \pi] $.
The initial density and velocity appear on figure \ref{ct5b_0_dv} and \ref{ct5b_2_dv}, together 
 with the numerical solutions computed at time $ t = 1 $ by means of the EPB and REPB-based schemes.
Table \ref{error_ct5b_tab} shows the numerical error on the potential by comparison against a solution computed with the classical scheme
 on a fine grid.

\begin{table}[hbtp]
\begin{tabular}{|c|c|c|c|c|}
\hline
$ N $ & $ \varepsilon_{EPB}(\phi), \lambda = 1 $ & $ \varepsilon_{REPB}(\phi), \lambda = 1 $ &
        $ \varepsilon_{EPB}(\phi), \lambda = 10^{-2} $ & $ \varepsilon_{REPB}(\phi), \lambda = 10^{-2} $ \\
\hline
\hline
$8000$  & $ 5 \times  10^{-5} $&$ 1.4 \times 10^{-4} $&$ 4.0 \times 10^{-4} $&$ 7.6 \times 10^{-4} $ \\
$4000$  & $ 1.5 \times 10^{-4} $&$ 2.8 \times 10^{-4} $&$ 1.0 \times 10^{-3} $&$ 1.8 \times 10^{-3} $ \\
$2000$  & $ 3.4 \times 10^{-4} $&$ 5.4 \times 10^{-3} $&$ 2.3 \times 10^{-3} $&$ 3.7 \times 10^{-3}$  \\
\hline
\end{tabular}
\caption{\label{error_ct5b_tab} Five-branch test case : 
comparison of the error in $ L^{\infty} $ norm
 on the potential for various grid sizes and dimensionless Debye length $ \lambda $
 using the EPB and REPB-based schemes.}
\end{table}

\subsubsection{Seven-branch solution}
\label{subsubsec_test2}

\begin{figure}[hbtp]
 \begin{minipage}[c]{.46\linewidth}
 \includegraphics[scale=0.4]{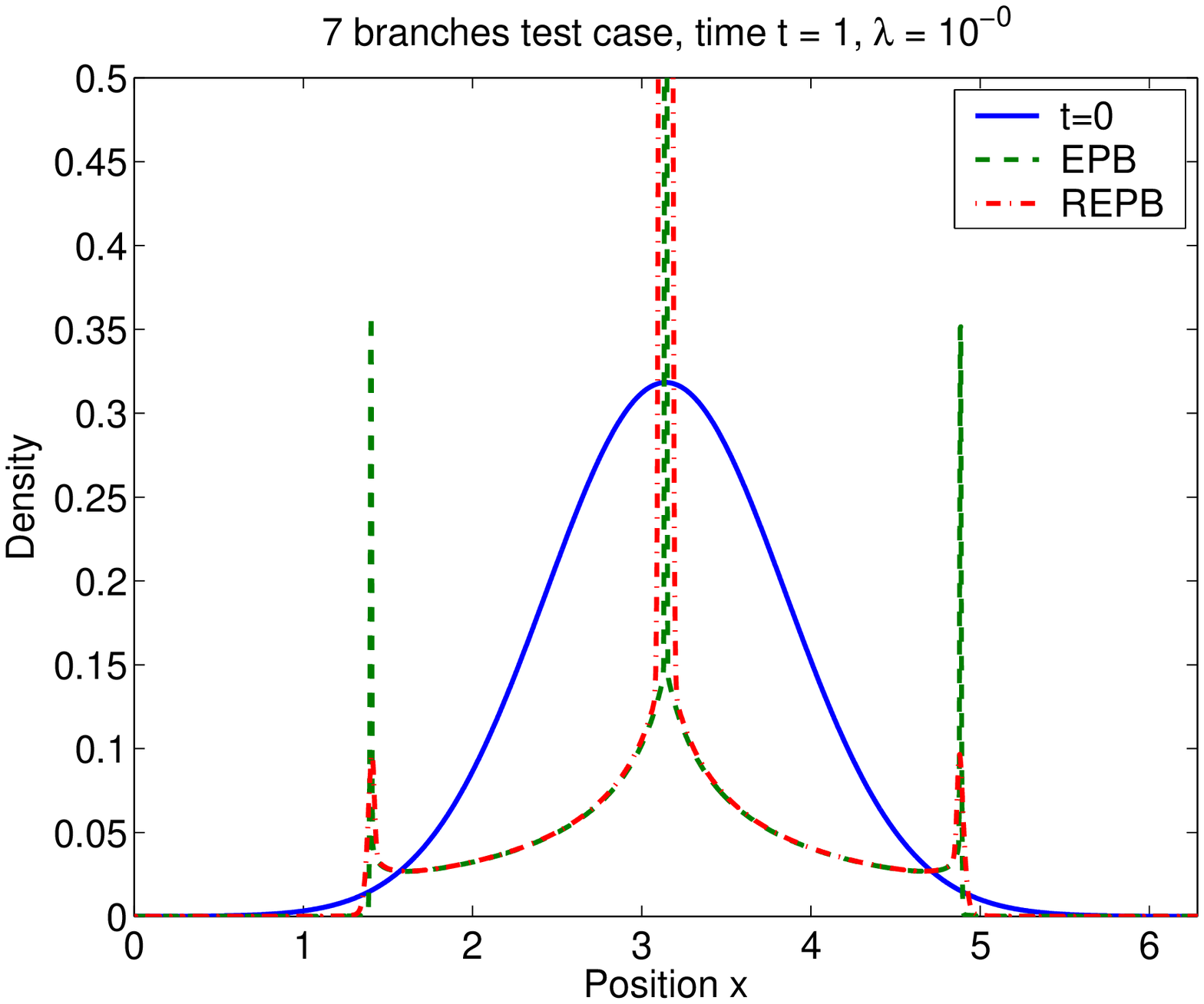}
 \end{minipage}
 \begin{minipage}[c]{.46\linewidth}
 \psfrag{varepsilon}{\footnotesize{$\varepsilon(nu)$}}
 \psfrag{Space step}{\footnotesize{$\Delta x$}}
 \includegraphics[scale=0.4]{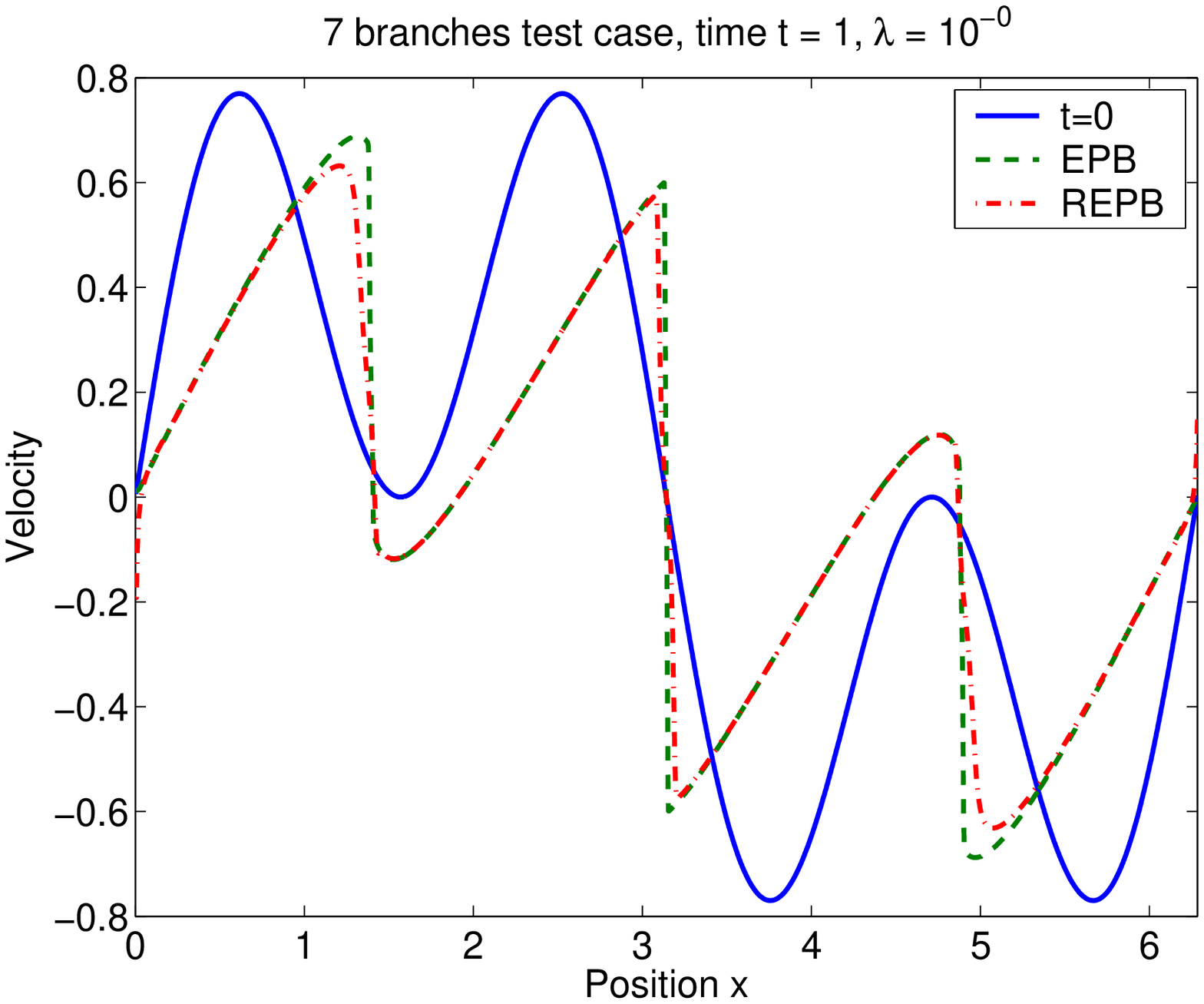}
 \end{minipage}
\caption{\label{ct7b_0_dv}Density (left) and velocity (right) as a function of space for the seven-branch test case.
The dimensionless Debye length $\lambda $ is $ 1 $. Numerical solutions at time $ t = 1 $ are computed with EPB and REPB-based
 schemes on a grid with $ 2000 $ cells.
}
 \begin{minipage}[c]{.46\linewidth}
 \psfrag{varepsilon}{\footnotesize{$\varepsilon(\phi)$}}
 \psfrag{Space step}{\footnotesize{$\Delta x$}}
 \includegraphics[scale=0.4]{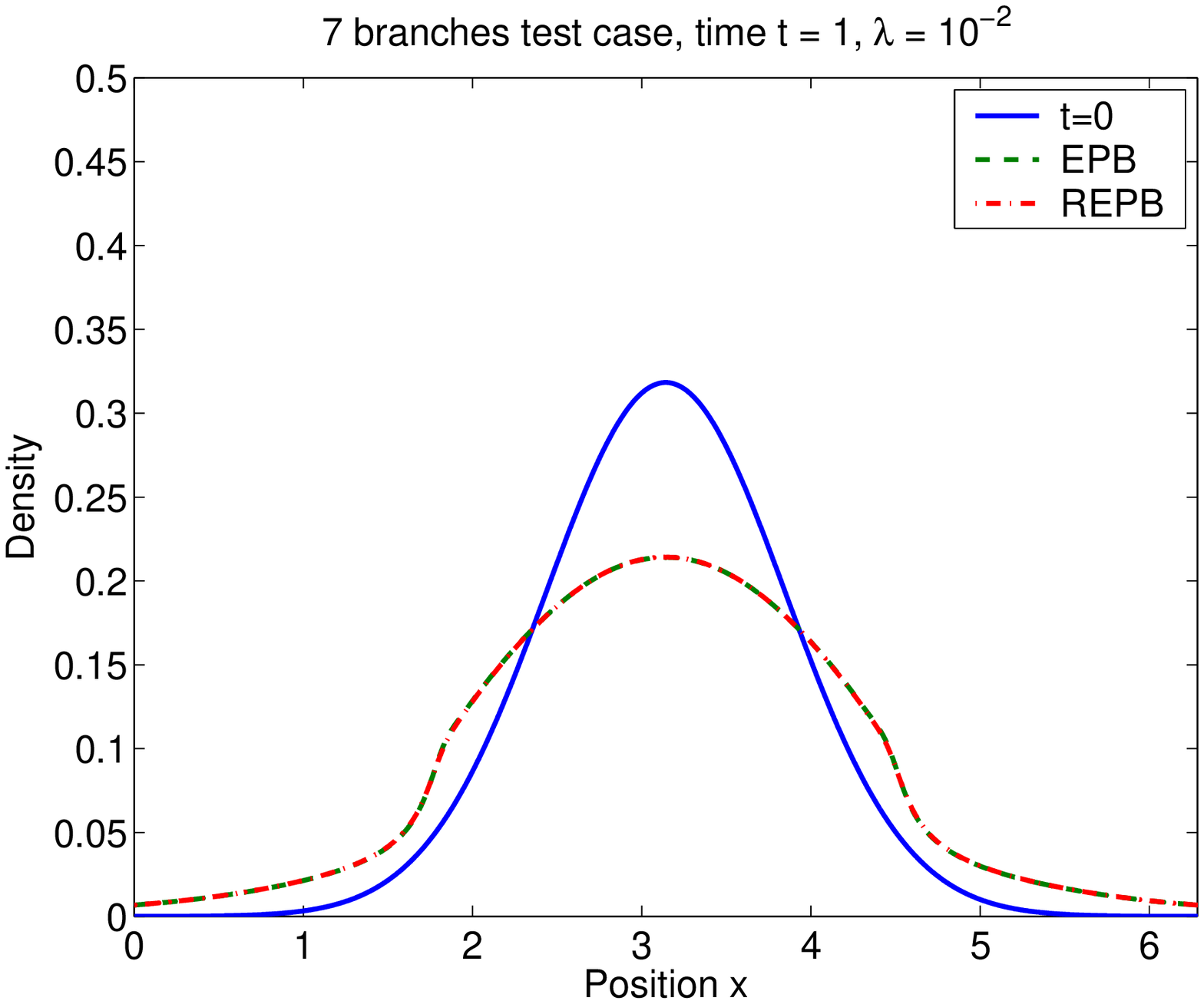}
 \end{minipage}
 \begin{minipage}[c]{.46\linewidth}
\includegraphics[scale=0.4]{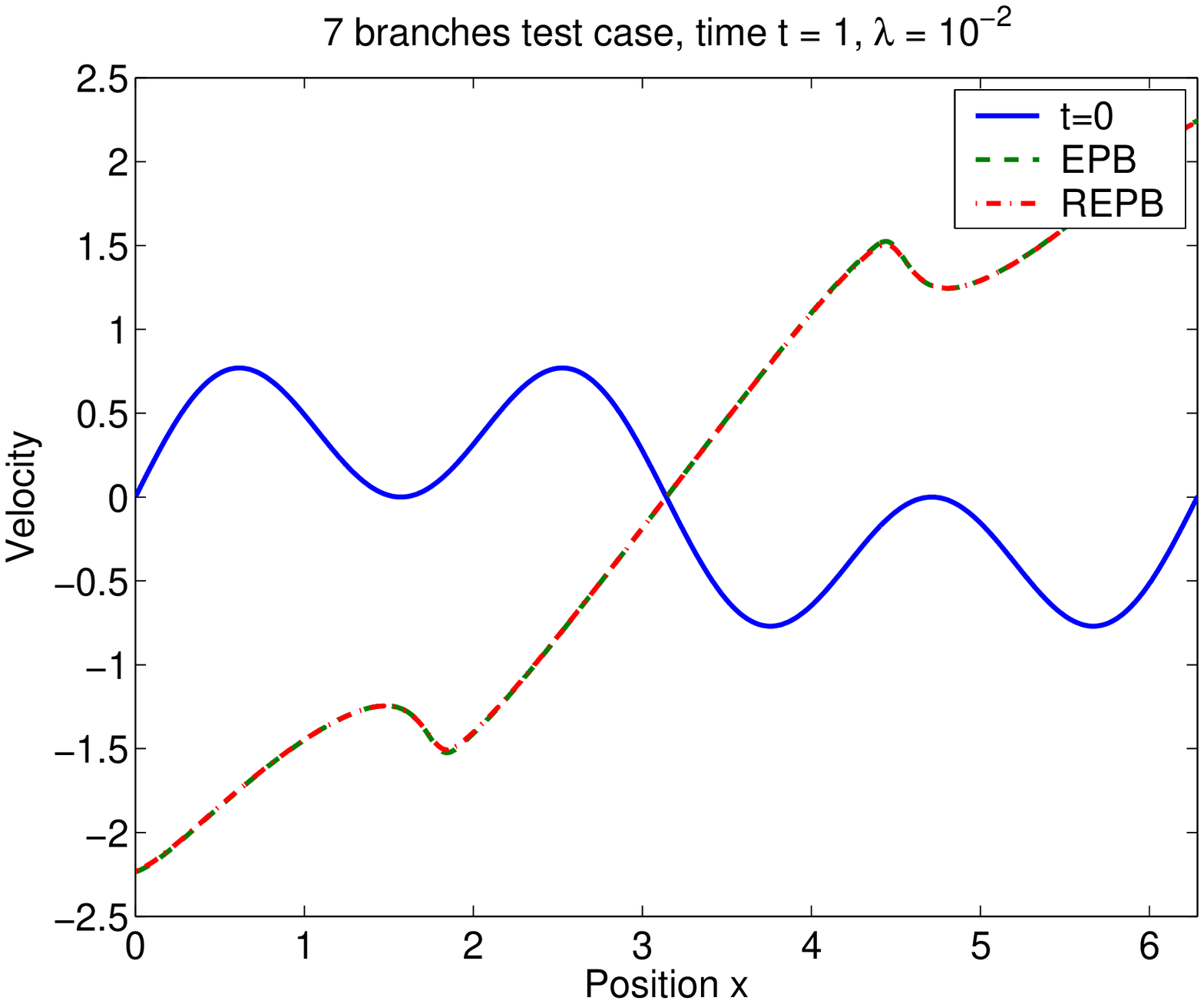}
 \end{minipage}
\caption{\label{ct7b_2_dv}Density (left) and velocity (right) as a function of space for the seven-branch test case.
The dimensionless Debye length $ \lambda $ is $ 10^{-2} $. Numerical solutions at time $ t = 1 $ are computed with EPB and REPB-based
 schemes on a grid with $ 2000 $ cells.}
\end{figure}

The initial condition for this second test-case is given by
\begin{eqnarray*}
& & n_{0}  = \frac{1}{\pi} e^{ -(x-\pi)^2 }, \\
& & u_{0}  = \sin(2x) \cos x.
\end{eqnarray*}
The initial density and velocity appear on figure \ref{ct7b_0_dv} and \ref{ct7b_2_dv}, together
 with the numerical solutions computed at time $ t = 1 $ using the EPB and REPB-based schemes.
Table \ref{error_ct7b_tab} shows the numerical error on the potential by comparison against a solution computed with the classical scheme
 on a fine grid.

\begin{table}[hbtp]
\begin{tabular}{|c|c|c|c|c|}
\hline
$ N $ & $ \varepsilon_{EPB}(\phi), \lambda = 1 $ & $ \varepsilon_{REPB}(\phi), \lambda = 1 $ &
        $ \varepsilon_{EPB}(\phi), \lambda = 10^{-2} $ & $ \varepsilon_{REPB}(\phi), \lambda = 10^{-2} $ \\
\hline
\hline
$8000$  & $ 4.6 \times  10^{-5} $&$ 4.5 \times 10^{-4} $&$ 2.3 \times 10^{-4} $&$ 7.6 \times 10^{-4} $ \\
$4000$  & $ 1.4 \times 10^{-4} $&$ 7.4 \times 10^{-4} $&$ 6.7 \times 10^{-4} $&$ 1.6 \times 10^{-3} $ \\
$2000$  & $ 3.2 \times 10^{-4} $&$ 1.2 \times 10^{-3} $&$ 1.29 \times 10^{-3} $&$ 3.5 \times 10^{-3}$  \\
\hline
\end{tabular}
\caption{\label{error_ct7b_tab} Seven-branch test case : comparison of the error in $ L^{\infty} $ norm
 on the potential for various grid sizes and dimensionless Debye length $ \lambda $
 using the classical and reformulated Euler-Poisson-Boltzmann schemes.}
\end{table}

\subsubsection{Analysis of the results for the five and seven branch test cases}
\label{subsubsec_dispersive_analysis}

In both the dispersive ($\lambda = 1 $) and hydrodynamic ($\lambda = 10^{-2}$) regimes, the two schemes give similar results.
On figure \ref{ct5b_0_dv}, \ref{ct5b_2_dv}, \ref{ct7b_0_dv} and \ref{ct7b_2_dv}, overlapping lines for the solution at time $ t = 1 $
 computed with the EPB and REPB-based schemes confirm this similar behavior.
One can exhibit some differences thanks to the numerical convergence study.
Tables \ref{error_ct5b_tab} and \ref{error_ct7b_tab} show the same differences between the two schemes as in the previously discussed soliton test case.
The error of the reformulated scheme is slightly larger than that of the classical scheme,
 and this difference is more obvious when $ \lambda = 1 $.


\setcounter{equation}{0}
\section{Conclusion}
\label{sec_conclu}

In this paper, we have analyzed two schemes for the Euler-Poisson-Boltzmann (EPB) model of plasma physics, and compared them in different regimes characterized by different values of the dimensionless Debye length $\lambda$. The dispersive regime corresponds to $\lambda = O(1)$ while the hydrodynamic regime is characterized by $\lambda \ll 1$. When $\lambda \to 0$, the EPB model formally converges to the Isothermal Compressible Euler (ICE) model. The first scheme we have considered is based on the original EPB formulation of the model. The second one uses a reformulation (referred to as the REPB model) in which the model more explicitly appears as a singular perturbation of the ICE Model. 

We have provided a stability analysis of the two schemes, showing that both schemes are stable in both the dispersive and hydrodynamic regimes, with stability constraints on the time and mesh steps which are independent of $\lambda$ when $\lambda \to 0$. Finally, we have tested them on three different one-dimensional test problems. The first test problem, the soliton test, provides an analytical solution in the dispersive regime. The second test problem, the Riemann problem with two expanding shock waves, is suitable to explore the hydrodynamic regime. Finally, the third test problem involves singularity formation in the dispersive regime. 

We have concluded that both scheme have similar behavior in the dispersive regime (with a slightly increased, but perfectly acceptable numerical diffusion in the case of the REPB-based schemes). By contrast, in the hydrodynamic regime, the EPB-based schemes develop oscillations and singularities, which, in under-resolved situations (i.e. when the time and space steps are too large to resolve the spatio-temporal variations of the solution) prevent any grid convergence of the solution. By contrast, the REPB-based scheme well captures the entropic solution of the ICE model, which provides a good approximation of the weak limit of the dispersive (oscillatory) solutions of the EPB model in the small $\lambda$ regime. 

Future works concern the extension of this analysis to the two- or multi-dimensional case, the passage to second order schemes and the pursuit of the analytical investigations of the accuracy and stability of the schemes in both regimes.


\end{document}